\documentclass{article}
\usepackage{amsmath}[1996/11/01]
\usepackage{amssymb,amsthm,amsxtra}
\usepackage{epsfig}

\def\[#1\]{\begin{equation}#1\end{equation}}
\makeatletter
\def\beq{%
   \relax\ifmmode
      \@badmath
   \else
      \ifvmode
         \nointerlineskip
         \makebox[.6\linewidth]%
      \fi
      $$%%$$ BRACE MATCH HACK
   \fi
}
\def\eeq{%
   \relax\ifmmode
      \ifinner
         \@badmath
      \else
         $$%%$$ BRACE MATCH HACK
      \fi
   \else
      \@badmath
   \fi
   \ignorespaces
}

\def\enddisplaymath{\eeq\global\@ignoretrue}
\makeatother

\newtheorem{thm}{Theorem}
\newtheorem{cor}[thm]{Corollary}
\newtheorem{lem}[thm]{Lemma}
\newtheorem{prop}[thm]{Proposition}
\newtheorem{rem}[thm]{\it{Remark}}
\newtheorem{rems}[thm]{\it{Remark}}

\numberwithin{equation}{section}
\numberwithin{thm}{section}
\DeclareMathOperator{\Exp}{\mathbb E}

\DeclareMathOperator{\Prob}{\mathbb P}

\DeclareMathOperator{\Plan}{Plan}
\DeclareMathOperator{\GUE}{GUE}
\DeclareMathOperator{\Pois}{Pois}

\newcommand{\C}{\mathbb C}
\newcommand{\R}{\mathbb R}
\newcommand{\Z}{\mathbb Z}
\newcommand{\N}{\mathbb N}

\newcommand{\calO}{\mathcal O}

\DeclareMathOperator{\K_n}{K_n}
\DeclareMathOperator{\Kk}{K}
\DeclareMathOperator{\Km_n}{K^{(m)}_n}

\DeclareMathOperator{\KN1}{K_{N+1}}

\DeclareMathOperator{\tK_N}{\tilde{K}_{N+1}}
\DeclareMathOperator{\ttK}{\tilde{K}_0}
\DeclareMathOperator{\S_n}{S_n}
\DeclareMathOperator{\Ss}{S}
\DeclareMathOperator{\Ssmn}{S^{(m)}_n}
\DeclareMathOperator{\Rrmn}{R^{(m)}_n}
\DeclareMathOperator{\Ssm}{S^{(m)}}
\DeclareMathOperator{\RR_n}{R_n}
\DeclareMathOperator{\Rr}{R}
\DeclareMathOperator{\Rrm}{R^{(m)}}
\DeclareMathOperator{\Tt}{T_n}

\begin{document}

\title{{\bf A Fredholm Determinant Identity
and \\ the Convergence of Moments \\ for Random Young Tableaux}}
\author{{\bf Jinho Baik,}\footnote{
Department of Mathematics,
Princeton University, Princeton, New Jersey, 08544,
jbaik@math.princeton.edu}
\footnote{Institute for Advanced Study,
Princeton, New Jersey 08540}
 \ \ {\bf Percy Deift}
\footnote{University of Pennsylvania, Philadelphia, 19104,
deift@math.upenn.edu}
\footnote{Courant Institute of Mathematical Sciences, New York 10012}
 \ \
and \ \ {\bf Eric M. Rains}\footnote{AT\&T Research, New Jersey,
Florham Park, New Jersey 07932, rains@research.att.com}}

\maketitle

\begin{abstract}

 We obtain an identity between Fredholm determinants 
of two kinds of operators, one acting on functions on the unit circle 
and the other acting on functions on a subset of the integers.
This identity is a generalization of an identity between 
a Toeplitz determinant and a Fredholm determinant that has appeared 
in the random permutation context.
Using this identity, we prove, in particular, convergence of moments 
for arbitrary rows of a random Young diagram under Plancherel measure.

\end{abstract}

%\tableofcontents

%%  revised 01/08/2001

\section{Introduction}\label{sec-intro}

In \cite{BDJ}, the authors considered the length $\ell_N(\pi)$ of the
longest increasing subsequence of a random permutation $\pi\in S_N$,
the symmetric group on $N$ numbers.
They showed, in particular, that
for $\tilde{\ell}_N(\pi):= \frac{\ell_N(\pi)-2\sqrt{N}}{N^{1/6}}$,
\begin{equation}\label{ei1}
  \lim_{N\to\infty} \Prob(\tilde{\ell}_N\le x) = F^{(1)}(x),
\end{equation}
where $F^{(1)}(x)$ is the Tracy-Widom distribution \cite{TW1} for the
largest eigenvalue of a random matrix from
the Gaussian Unitary Ensemble (GUE).
The authors also proved the convergence of moments,
\begin{equation}\label{ei2}
  \lim_{N\to\infty} \Exp\bigl((\tilde{\ell}_N)^m\bigr)
= \int_{-\infty}^\infty x^m d F^{(1)}(x),
\qquad m=1,2,\cdots.
\end{equation}
The authors then reinterpreted \eqref{ei1}, \eqref{ei2} in terms of
Young diagrams $\lambda=(\lambda_1,\lambda_2,\cdots)$
via the Robinson-Schensted correspondence.
Here $\lambda_j$ is the number of boxes in the $j$th row of $\lambda$
and $\lambda_1\ge \lambda_2\ge\cdots\ge 0$.
%(A reference for Young diagrams and Young tableaux is
%\cite{Stl}.)
The set of Young diagrams $Y_N$ of size $N$, $\sum_j \lambda_j =N$,
is equipped with \emph{Plancherel measure},
\begin{equation}\label{ei3}
  \Prob^{\Plan}_N(\lambda) := \frac{d_\lambda^2}{N!}, \qquad \lambda\in Y_N,
\end{equation}
where $d_\lambda$ is the number of standard Young tableaux of shape $\lambda$.
Set
\begin{equation}\label{ei4}
   \xi_j:= \frac{\lambda_j-2\sqrt{N}}{N^{1/6}}, \qquad j=1,2,\cdots.
\end{equation}
Then \eqref{ei1}, \eqref{ei2} imply that $\xi_1$ converges in distribution,
together with all its moments, to $F^{(1)}$.
This reinterpretation led the authors to conjecture that
for all $k$, $\xi_1, \xi_2, \cdots, \xi_k$ converge
to the joint distribution function $F(x_1,x_2, \cdots, x_k)$
for the first $k$ eigenvalues of a random GUE matrix.
%taken from the Gaussian unitary ensemble (GUE).
In \cite{BDJ2}, the authors verified the convergence in distribution,
together with its moments, to the Tracy-Widom distribution
$F^{(2)}$ for the second largest eigenvalue of a random GUE matrix.
The conjecture for $\xi_1,\xi_2,\cdots, \xi_k$ was then proved
in three independent papers \cite{Ok}, \cite{BOO}, \cite{kurtj:disc},
all appearing within a few months in the spring of 1999.
Let $y_j$ be the $j$th largest eigenvalue of
a random $N\times N$ matrix from GUE with probability density
\begin{equation}
   d\Prob^{\GUE}_N(y_1,\cdots, y_N)
=\frac1{Z_N} \prod_{1\le i<j\le N} (y_i-y_j)^2 \prod_{j=1}^N e^{-y_j^2}
dy_1\cdots dy_N,
\end{equation}
where $y_1\ge \cdots\ge y_N$, and $Z_N$ is the normalization constant.
At the `edge' of the spectrum, the following convergence in distribution
is well-known (see, e.g. \cite{TW1}, \cite{kurtj:disc} Theorem 1.4):
for any $k\in\N$,
there is a distribution function $F(x_1,\cdots,x_k)$
on $x_1\ge\cdots\ge x_k$ such that
\begin{equation}\label{eq2.4}
  \lim_{N\to\infty}\Prob^{\GUE}_N\bigl((y_1-\sqrt{2N})\sqrt{2}N^{1/6}\le x_1,
\cdots, (y_k-\sqrt{2N})\sqrt{2}N^{1/6}\le x_k\bigr)
= F(x_1,\cdots,x_k).
\end{equation}
In all three papers \cite{Ok}, \cite{BOO}, \cite{kurtj:disc},
the authors showed that for any $x_1,\cdots, x_k\in \R^k$,
\begin{equation}\label{ei5}
   \lim_{N\to\infty} \Prob^{Plan}_N(\xi_1\le x_1, \cdots, \xi_k\le x_k)
= F(x_1,\cdots,x_k),
\end{equation}
but the question of the convergence of moments was left open.

Introduce the \emph{Poissonized Plancherel measure}
\begin{equation}
  \Prob^{Pois}_t(\lambda) = \sum_{N=0}^\infty
\frac{e^{-t^2}t^{2N}}{N!}\Prob^{Plan}_N(\lambda),
\qquad t>0,
\end{equation}
on all Young diagrams, which corresponds to choosing
$N$ as a Poisson variable with parameter $t^2$.
Here $\Prob^{\Plan}_N(\lambda)=0$ if $\lambda$ is not a partition of $N$.
Throughout the paper, we will work with $\Prob^{Pois}_t(\lambda)$
rather than $\Prob^{\Plan}_N(\lambda)$ itself.
This is because
the expectation with respect to $\Prob^{\Pois}_t(\lambda)$
leads to convenient determinantal formulae.
Indeed, in \cite{Gessel}, Gessel proved the following formula
\begin{equation}\label{ei7}
   \Prob^{Pois}_t(\lambda_1\le n)= e^{-t^2} \det(\Tt),
\end{equation}
where $\Tt$ is the $n\times n$ Toeplitz matrix
with entries $(\Tt)_{pq}=c_{p-q}$, $0\le p,q <n$,
where $c_k$ is the $k^{th}$ Fourier coefficient
of $e^{t(z+z^{-1})}$,
$c_k=\int_{|z|=1} z^{-k}e^{t(z+z^{-1})}\frac{dz}{2\pi iz}$.
This formula played a basic role in \cite{BDJ} in proving
\eqref{ei1}, \eqref{ei2}.
In \cite{BDJ2}, the authors introduced the integral operator $\K_n$
with $\varphi(z)=e^{t(z-z^{-1})}$ (see \eqref{3} below)
and proved the following formulae
\begin{equation}\label{ei8}
  \Prob^{Pois}_t(\lambda_1\le n)= 2^{-n}\det(1-\K_n)
\end{equation}
and
\begin{equation}\label{ei9}
 \Prob^{Pois}_t(\lambda_2\le n+1)=
\Prob^{Pois}_t(\lambda_1\le n)
+ \biggl(-\frac{\partial}{\partial s}\biggr)\bigg|_{s=1}
[(1+\sqrt{s})^{-n}\det(1-\sqrt{s}\K_n)].
\end{equation}
These formulae played a basic role in \cite{BDJ2} in proving
the analogue of \eqref{ei1},\eqref{ei2} for $\lambda_2$.
In \cite{BOO} and \cite{kurtj:disc}, and also later, in greater
generality, in \cite{Ok99} and \cite{Rains:corr},
the authors obtained the following identity :
Let $\Lambda_k$ denote the (finite) set
$\{n\in\{0,1,\cdots\}^k : \sum_{j=1}^r n_j\le r-1, r=1,\cdots, k\}$.
Then for $a_k\le \cdots\le a_1\le a_0=\infty$,
\begin{equation}\label{ei10}
\begin{split}
  &\Prob^{Pois}_t(\lambda_1-1\le a_1, \lambda_2-2\le a_2, \cdots,
\lambda_k-k\le a_k) \\
&\qquad =
\sum_{n\in\Lambda_k} \frac1{n_1!\cdots n_k!}
\frac{\partial^{|n|}}{\partial s_1^{n_1}\cdots \partial s_k^{n_k}}
\bigg|_{s_1=\cdots=s_k=-1}
\det(1+(\sum_{l=1}^k s_l\chi_{(a_l,a_{l-1}]})\Ss),
\end{split}
\end{equation}
where the matrix elements of $\Ss(i,j)$ are given in \eqref{ef1} below
with $\varphi(z)=e^{\sqrt{\gamma}(z-z^{-1})}$.
As usual, $\chi_{(a,b]}$ denotes the characteristic function
of the interval $(a,b]$, and so
%Here $\Ss$ has the matrix elements
%\begin{equation}
%  \Ss(i,j)=\sum_{k\ge 1} (\varphi^{-1})_{i+k}\varphi_{-j-k},
%\qquad i,j\in\Z,
%\end{equation}
%where $\varphi(z)=e^{\sqrt{\gamma}(z-z^{-1})}=\sum_{j\in\Z}\varphi_j z^j$,
%$\varphi^{-1}(z)=e^{-\sqrt{\gamma}(z-z^{-1})}
%=\sum_{j\in\Z}(\varphi^{-1})_j z^j$.
%Thus
$(\sum_{l=1}^k s_l\chi_{(a_l,a_{l-1}]})\Ss$ denotes the operator
in $\ell^2(\Z)$ with kernel $s_l\Ss(i,j)$ if $i\in(a_l,a_{l-1}]$,
and zero otherwise.
Setting $a_j=2t+x_jt^{1/3}$, $x_1\ge x_2\ge\cdots\ge x_k$,
and letting $t\to\infty$, and de-Poissonizing as in \cite{Jo2},
the authors in \cite{BOO} and \cite{kurtj:disc} obtain \eqref{ei5}.
In \cite{BOO} and \cite{kurtj:disc}, however, the authors
are not able to prove convergence of moments.
The reason for this is that it is possible to use the classical
steepest-descent method to control
$\det(1+(\sum_{l=1}^k s_l\chi_{(a_l,a_{l-1}]})\Ss)$
for $a_j=2t+x_jt^{1/3}$ as $t\to\infty$,
uniformly for $x_1\ge x_2\ge\cdots\ge x_k\ge M$ for any
fixed $M$.
But as the $x_j$'s tends to $-\infty$,
the method break down.
On the other hand,
the authors in \cite{BDJ, BDJ2} are able to control the lower tails
of the probability distributions,
and hence prove the convergence of moments for $\lambda_1$ and $\lambda_2$,
using the steepest-descent method for the Riemann-Hilbert problem
(RHP) naturally associated with $\Tt$ and $\K_n$ above.
The steepest-descent method for RHP was introduced in \cite{DZ1},
and extended to include fully non-linear oscillations in \cite{DVZ}.
The asymptotic analysis in \cite{BDJ}, \cite{BDJ2} is closely
related to the analysis in \cite{DKMVZ2, DKMVZ3}.
The main motivation for this paper was to find a formula for the joint
distribution of $\lambda_1,\cdots,\lambda_k$,
which generalized \eqref{ei9},
and to which the above Riemann-Hilbert steepest-descent methods
could be applied to obtain the lower tail estimates.

Note that from \eqref{ei7}, \eqref{ei8} and \eqref{ei10},
we have three formulae for the distribution of $\lambda_1$,
\begin{equation}\label{ei12}
\begin{split}
  \Prob^{Pois}_t(\lambda_1\le n) &= e^{-t^2}\det(\Tt)\\
&= 2^{-n} \det(1-\K_n)\\
&= \det(1-\chi_{[n,\infty)}\Ss)
\end{split}
\end{equation}
and from \eqref{ei9} and \eqref{ei10},
two formulae for the distribution of $\lambda_2$,
\begin{equation}\label{ei13}
\begin{split}
  \Prob^{Pois}_t(\lambda_2\le n+1)
& =\Prob^{Pois}_t(\lambda_1\le n)
+ \biggl(-\frac{\partial}{\partial s}\biggr)\bigg|_{s=1}
[(1+\sqrt{s})^{-n}\det(1-\sqrt{s}\K_n)]\\
&  =\Prob^{Pois}_t(\lambda_1\le n)
+ \biggl(\frac{\partial}{\partial s}\biggr)\bigg|_{s=-1}
\det(1+s\chi_{[n,\infty)}\Ss).
\end{split}
\end{equation}
To obtain the second formula, we use the fact that
$\Lambda_{k=2}=\{(0,0), (0,1)\}$ and set $a_1=\infty, a_2=n-1$
in \eqref{ei10}.
From \eqref{ei13}, we might guess that
\begin{equation}\label{ei14}
  (1+\sqrt{s})^{-n}\det(1-\sqrt{s}\K_n)=\det(1-s\chi_{[n,\infty)}\Ss).
\end{equation}
The content of Theorem \ref{mainthm} is that precisely this
relation is true for a general class of functions $\varphi(z)$,
provided $\varphi(z)$ has no winding. If the winding number of
$\varphi$ is non-zero, the above relation must be modified
slightly as in \eqref{ethm2}. The fact that
$e^{-t^2}\det(\Tt)=\det(1-\chi_{[n,\infty)}\Ss)$ for (essentially)
the same general class of $\varphi's$ (with zero winding number)
was first proved in \cite{BoOk}, with an alternative proof given
in \cite{BW}. The relation \eqref{ei14} for general $s$ was proved
essentially simultaneously with the present paper by Rains in
\cite{Rains:corr}, for a subclass of functions $\varphi$ with zero
winding, using algebraic methods (see Remark 4 in Section
\ref{sec-iden}). A particularly simple proof of the relation
$e^{-t^2}\det(\Tt)=\det(1-\chi_{[n,\infty)}\Ss)$ can be found in
the recent paper \cite{Bott2} of B\"ottcher (see also
\cite{Bott1}). The paper \cite{Bott2} also extends Theorem
\ref{mainthm} and \ref{thmmulti} to the matrix case (see Remark
\ref{rem-Bott1} and \ref{rem-Bott2} below).

In this paper, we will prove a general identity between determinants of
operators of two types :
the operators of the first type act on functions on the unit circle,
and the operator of the second type act on functions on a subset of
the integers.
Specializations of this identity have, in particular, the following
consequences :
\begin{enumerate}
\item[(S1)] A proof of the convergence of moments for $\xi_1,\cdots,\xi_k$
(see Theorem \ref{thmmom})
\item[(S2)] An interpretation of $F(x_1,\cdots, x_k)$ in \eqref{ei5}
as a ``multi-Painlev\`e'' function (see Section \ref{sec-ptons}).
As we will see, the behavior of multi-Painlev\'e functions has similarities
to the interactions of solitons in the classical theory of
the Korteweg de Vries equation.
\item[(S3)] The analogue of Theorem \ref{thmmom} for signed permutations
and so-called colored permutations (see Section \ref{sec-color})
\item[(S4)] New formulae for random word problems, certain 2-dimensional growth
models,  and also the so-called ``digital boiling'' model
(see Section \ref{sec-color})
\end{enumerate}
The new identity is given in Theorem \ref{mainthm} in two closely related
forms \eqref{ethm2}, \eqref{ethm2.1}.
In (S1)-(S4), we only use \eqref{ethm2}.

As we will see, some simple estimates together
with a Riemann-Hilbert analysis of $\det(1-\sqrt{s}\K_n)$
is enough to control the lower tail estimation of
$\Prob^{\Pois}_t(\lambda)$.
The relation \eqref{ei14} generalizes to the multi-interval case, as
described in Theorem \ref{thmmulti} in Section \ref{sec-iden}.
%This multi-interval relation leads to a ``multi-Painlev\`e''
%interpretation
%of $F(x_1,\cdots, x_k)$ as described in Section \ref{sec-ptons}.

In Section \ref{sec-iden}, we prove the main identity \eqref{ethm2},
\eqref{ethm2.1} in the single interval case, and also
the identity \eqref{ethmmulti} in the multi-interval case.
In Section \ref{sec-mom}, we use \eqref{ethm2} to prove the convergence
of moments for random Young tableaux (Theorem \ref{thmmom}).
A stronger version of this result is given in \eqref{ewq}.
Section \ref{sec-tail} contains certain tail estimates, needed
in Section \ref{sec-mom}.
Various estimates needed in Section \ref{sec-tail} for a ratio
of determinants are derived in Section \ref{sec-rhp} using the
steepest-descent method for RHP's.
In Section \ref{sec-ptons}, we introduce the notion of a
multi-Painlev\'e solution, and in Section \ref{sec-color}, we prove
various formulae for colored permutations and also discuss certain
random growth models from the perspective of Theorem \ref{mainthm}.

\medskip
\noindent {\bf Acknowledgments.} The authors would like to thank
Xin Zhou for useful comments. The authors would also like to thank
Albrecht B\"ottcher for pointing out a calculational error in an
earlier version of the text. The work of the first author was
supported in part by NSF Grant \# DMS 97-29992. The work of the
second author was supported in part by NSF Grant \# DMS 00-03268,
and also by the Guggenheim Foundation.

%%  revised 01/08/2001

\section{Fredholm determinant identity}\label{sec-iden}

%Consider a function $\psi(z)$ on the unit circle $\Sigma=\{ z\in\C : |z|=1\}$.
%We assume that it can be written as
%\begin{equation}\label{1}
%  \psi(z)=\psi_+(z)\psi_-(z)
%\end{equation}
%with $\psi_+(0)=\psi_-(\infty)=1$ and $\psi_+$ (resp. $\psi_-$)
%extends to a nonzero analytic function in the interior (resp. exterior)
%of the circle.
%We define
%\begin{equation}\label{2}
%  \varphi(z)= \psi_+(z)\psi_-(z)^{-1},
%\end{equation}

Let $\varphi(z)$ be a continuous, complex-valued, non-zero
function on the
unit circle $\Sigma=\{ z\in\C : |z|=1\}$.
Define $\K_n$ to be the integral operator acting on $L^2(\Sigma,dw)$
with kernel
\begin{equation}\label{3}
  \K_n(z,w)=\frac{1-z^n\varphi(z)w^{-n}\varphi(w)^{-1}}{2\pi i(z-w)},
\qquad (\K_nf)(z)=\int_{|w|=1} \K_n(z,w)f(w)dw.
\end{equation}
For a function $f$ on $\Sigma$, its Fourier coefficients are denoted by
$f_j$, so that
\begin{equation}
  f(z)=\sum_{j\in\Z} f_j z^j.
\end{equation}
Let $\Ss$ be the matrix with entries
\begin{equation}\label{ef1}
  \Ss(i,j)=\sum_{k\ge 1} (\varphi^{-1})_{i+k}\varphi_{-j-k},
\qquad i,j\in\Z,
\end{equation}
and let $\Rr$ be the matrix with entries
\begin{equation}\label{ef2}
  \Rr(i,j)=\sum_{k\le 0} (\varphi^{-1})_{i+k}\varphi_{-j-k},
\qquad i,j\in\Z,
\end{equation}
Let $\S_n$ denote the operator $\chi_{[n,\infty)}\Ss$
acting on $\ell^2(\{n,n+1,\cdots\})$,
%Define the operator $\S_n$ acting on $\ell^2(\{n,n+1,\cdots\})$
%with matrix elements
\begin{equation}\label{en0}
%  \S_n(i,j)=\sum_{k\ge 1} (\varphi^{-1})_{i+k}\varphi_{-j-k},
%\qquad
(\S_nf)(i)=\sum_{j\ge n} \S_n(i,j)f(j), \qquad i\ge n,
\end{equation}
and let $\RR_n$
denote the operator $\chi_{(-\infty, n-1]}\Rr$
acting on $\ell^2(\{\cdots,n-2,n-1\})$,
%with matrix elements
\begin{equation}%\label{en0}
%  \RR_n(i,j)=\sum_{k\le 0} (\varphi^{-1})_{i+k}\varphi_{-j-k},
%\qquad
(\RR_nf)(i)=\sum_{j\le n-1} \RR_n(i,j)f(j), \qquad i\le n-1.
\end{equation}

\begin{thm}\label{mainthm}
Let $\varphi(z)$ be a non-zero function on the unit circle
satisfying $\sum_{j\in\Z} |j\varphi_j| <\infty$,
which has winding number equal to $\#(\varphi)$.
For $s\in\C$ and $n\in\Z$,
$\K_n$, $\S_n$ and $\RR_n$ are trace class on $L^2(\Sigma,dw)$,
$\ell^2(\{n,n+1,\cdots\})$ and $\ell^2(\{\cdots, n-1\})$ respectively, and
we have %for $s\neq -1$,
\begin{eqnarray}\label{ethm2}
  \det (1-s\K_n) &=& (1+s)^{n+\#(\varphi)}\det (1- s^2\S_n),
\qquad s\neq -1,\\ \label{ethm2.1} &=& (1-s)^{-n-\#(\varphi)}\det
(1- s^2\RR_n), \quad s\neq 1.
\end{eqnarray}
\end{thm}

\begin{rems}
Standard Banach algebra estimates show that
if the winding number of $\varphi$ is zero and
$\sum |j\varphi_j|<\infty$,
then $\|\log\varphi\|_{\infty} + (\sum|j||(\log\varphi)_j|^2)^{1/2}<\infty$.
This is enough to prove that the first and the third terms in \eqref{ei12}
are equal for all such $\varphi$'s (see \cite{BW}).
In particular, by \eqref{ei12}, \eqref{ethm2} is true
for all $\varphi$ without winding and satisfying $\sum |j\varphi_j|<\infty$,
when $s=1$.
\end{rems}

\begin{rems}\label{rem-Bott1}
  As noted by B\"ottcher \cite{Bott2}, Theorem \ref{mainthm}
  remains true in the case where $\varphi(z)$ is an invertible
  $N\times N$ matrix, provided the exponent $n+\#(\varphi)$ is
  replaced in \eqref{ethm2}, \eqref{ethm2.1} by
  $Nn+\#(\det\varphi)$. The proof in the scalar case extends to
  $N\times N$ matrices, and we give no further details : the proof
  in \cite{Bott2} is different and uses Wiener-Hopf factorization
  directly.
\end{rems}

For the proof of Theorem \ref{mainthm},
we use the following basic properties of
the determinant (see, e.g., \cite{SimonTr}).
If $A$ is a trace class operator on a Hilbert space $H$,
$\|A\|_1= (tr A^*A)^{1/2}$ denotes the trace norm.

\begin{lem}\label{lem0}
\begin{enumerate}
\item\label{lem1}
 If $A_n$ is a trace class operator for each $n$ and $A_n\to A$ in trace norm,
then $A$ is a trace class operator and $\det(1+A_n)\to \det(1+A)$
as $n\to\infty$.
\item\label{lem2}
 If $A$ is a trace class operator, and $B_n$ and $C_n$
are bounded operators such that $(B_n)^*$ and $C_n$ converge
strongly to $B^*$ and $C$ respectively, then $\det(1+C_nAB_n)\to
\det(1+CAB)$ as $n\to\infty$.
\item\label{lem3}
 If $AB$ and $BA$ are trace class operators,
then $\det(1+AB)=\det(1+BA)$.
\item\label{lem4}
 Suppose $C$ acts on $\ell^2(\Z)$ and has matrix elements
$(c_{ij})_{i,j\in\Z}$.
If $\sum_{i,j\in\Z} |c_{ij}|<\infty$, then $C$ is trace class
and $\|C\|_1\le \sum_{i,j\in\Z} |c_{ij}|$.
\end{enumerate}
\end{lem}

\begin{proof}[Proof of Theorem \ref{mainthm}]
Define the projection operators on the circle
\begin{equation}
  (P_nf)(z)= \sum_{j\ge n} f_jz^j, \qquad n\in\Z,
\end{equation}
and
\begin{equation}
\begin{split}
  &(Q_nf)(z)=\sum_{0\le j< n} f_jz^j, \quad n>0,\\
  &(Q_nf)(z)=-\sum_{n\le j< 0} f_jz^j, \quad n<0,\\
\end{split}
\end{equation}
with $(Q_0f)(z)=0$.
Thus in particular, we have $P_n=P_0-Q_n$.
Let $M_g$ denote the multiplication operator
\begin{equation}
  (M_gf)(z)=g(z)f(z).
\end{equation}
Direct calculation shows that
\begin{equation}\label{e9}
\begin{split}
  &\K_n = -P_0 + M_\varphi P_n M_{\varphi^{-1}}
=(1-P_0)-M_\varphi (1-P_n) M_{\varphi^{-1}}, \\
&\S_n=P_nM_{\varphi^{-1}}(1-P_0)M_\varphi P_n,\\
&\RR_n= (1-P_n)M_{\varphi^{-1}}P_0 M_\varphi (1-P_n).
%= -Q_n-P_n + M_\varphi P_n M_{\varphi^{-1}}.
\end{split}
\end{equation}

First, we show that $\K_n$, $\S_n$ and $\RR_n$ are trace class.
Indeed $\K_n=-Q_n-HM_{\varphi^{-1}}$,
where $H=[P_n, M_\varphi]$.
$H$ acts on the basis $\{z^l\}_{l\in\Z}$ for $L^2(\Sigma,dw)$,
as follows :
$Hz^k= \sum_{l} H_{lk} z^l$.
We find
\begin{equation}
  H_{lk}= \varphi_{l-k} (\chi_{l\ge n}- \chi_{k\ge n}), \qquad l,k\in\Z,
\end{equation}
where $\chi_{\cdot\ge n}$ denotes the characteristic function
of the set $\{ k\ge n\}$.
But $\sum_{l,k} |H_{lk}| \le \sum_{j} |j\varphi_j| <\infty$,
and hence by Lemma \ref{lem0} \ref{lem4},
we have the trace norm estimate
\begin{equation}\label{ep2}
  \|\K_n\|_1 \le n+
\biggl(\sum_{j} |j\varphi_j|\biggr) \| \varphi^{-1}\|_{L^\infty}.
\end{equation}
Now write $\S_n=AB$ where $A : \ell^2(\{1,2,\cdots \})
\to \ell^2(\{n,n+1,\cdots \})$
and $B : \ell^2(\{n,n+1,\cdots \})
\to \ell^2(\{1,2,\cdots \})$
with matrix elements
\begin{equation}
\begin{split}
 &A_{ik} =(\varphi^{-1})_{i+k},  \qquad i\ge n, k\ge 1, \\
 &B_{kj}=\varphi_{-k-j},  \quad\qquad k\ge 1, j\ge n.
\end{split}
\end{equation}
Write
\begin{equation}\label{ep4}
  A=  \chi^+_n \Phi^{-1} \chi^-_{-1} R
\end{equation}
where $(Rf)_j=f_{-j}$,
$\Phi^{-1}$ denotes convolution on $\ell^2(\Z)$ by $\{ \varphi^{-1}_j\}$,
\begin{equation}
  (\Phi^{-1} h)_j = \sum_{l\in\Z} (\varphi^{-1})_{j-l}h_l,
\end{equation}
and $\chi^+_n$, $\chi^-_{-1}$ are the projections
onto $\{k\ge n\}$ and $\{k\le -1\}$ respectively.
From \eqref{ep4}, it is clear that $A$ is bounded
from $\ell^2(\{1,2,\cdots \})
\to \ell^2(\{n,n+1,\cdots \})$
with norm estimate
\begin{equation}\label{ep5}
  \|A\| \le \|\varphi^{-1}\|_{L^\infty}.
\end{equation}
On the other hand, a similar calculation to \eqref{ep2} shows that
$B$ is trace class
from $\ell^2(\{n,n+1,\cdots \})
\to \ell^2(\{1,2,\cdots \})$ and
\begin{equation}\label{ep6}
  \| B\|_1 \le \sum_{l\ge n+1} |l\varphi_{-l}|
\le \sum_{l} |l\varphi_{l}|,
\end{equation}
which implies
\begin{equation}\label{ep7}
  \|\S_n\|_1 \le \biggl(\sum_{l\ge n+1} |(l+|n|)\varphi_{-l}|\biggr)
\|\varphi^{-1}\|_{L^\infty}.
\end{equation}
Similarly, we have
\begin{equation}\label{ep7.5}
  \|\RR_n\|_1 \le \biggl(\sum_{l\le n} |(l+|n|)\varphi_{-l}|\biggr)
\|\varphi^{-1}\|_{L^\infty}.
\end{equation}
Thus $\K_n$, $\S_n$ and $\RR_n$ are trace class.
Moreover, if we set $\varphi_J:=\sum_{|j|\le J} \varphi_jz^j$, $J\ge 0$,
then from (the proofs of) \eqref{ep2}, \eqref{ep7} and \eqref{ep7.5},
it is clear that as $J\to\infty$,
$\K_n(\varphi_J)\to \K_n(\varphi)$,
$\S_n(\varphi_J)\to \S_n(\varphi)$,
$\RR_n(\varphi_J)\to \RR_n(\varphi)$ in trace norm,
and hence the Fredholm determinants converge to
the corresponding determinants.
% $\det(1-s\K_n(\varphi_J)) \to
%\det(1-s\K_n(\varphi))$,
%$\det(1-s\S_n(\varphi_J)) \to \det(1-s\S_n(\varphi))$.
Also for $J$ sufficiently large, the winding number of $\varphi_J$
is the same as the winding number of $\varphi$,
and so we see that to prove \eqref{ethm2}, it is enough to consider
$\varphi$'s which are non-zero and analytic in a neighborhood of
$\Sigma$.
Henceforth we will assume that $\varphi$ is analytic :
this analyticity assumption
is not necessary and is used only to give a particularly
simple proof of Lemma \ref{bdj2} below.

In the below, we only present the proof of \eqref{ethm2}.
The proof of \eqref{ethm2.1} is similar.
Formally, we proceed as follows.
Suppose $P_n$ is finite rank so that $P_0=Q_n+P_n$
is also finite rank.
We have
\begin{equation}\label{e10}
  \det(1+s P_0 - s M_\varphi P_n M_{\varphi^{-1}})
= \det(1+s P_0)
\det(1-\frac{s}{1+sP_0}M_\varphi P_n M_{\varphi^{-1}})
\end{equation}
Using $P_0=Q_n+P_n$ and $\frac{1}{1+sP_0}=1-\frac{s}{1+s}P_0$,
the right-hand-side reduces to
\begin{equation}\label{e11}
\begin{cases}
  \det(1+s Q_n) \det(1+s P_n)
\det(1-s(1-\frac{s}{1+s}P_0)M_\varphi P_n M_{\varphi^{-1}}),
\quad &n\ge 0, \\
  \bigl( \det(1-s Q_n) \bigr)^{-1} \det(1+s P_n)
\det(1-s(1-\frac{s}{1+s}P_0)M_\varphi P_n M_{\varphi^{-1}}),
\quad &n< 0.
\end{cases}
\end{equation}
The first term in both cases is equal to $(1+s)^n$. Using Lemma
\ref{lem0} \ref{lem3} and $P_n=P_n^2$ for the last determinant,
\eqref{e11} becomes
\begin{equation}\label{e13}
\begin{split}
 &(1+s)^{n}\det(1+s P_n)
\det(1-sP_n M_{\varphi^{-1}}(1-\frac{s}{1+s}P_0)M_\varphi P_n )\\
&\quad =  (1+s)^{n}\det( (1+s P_n) - s(1+s) P_n
M_{\varphi^{-1}}(1-\frac{s}{1+s}P_0)M_\varphi P_n )\\ &\quad =
(1+s)^{n}\det( 1-s^2 \S_n),
\end{split}
\end{equation}
which is the desired result, up to the winding number $\#(\varphi)$.
For the case in hand, however, $P_0$ is not a trace class operator and the
above ``proof'' breaks down.
We circumvent the difficulty by approximating the operator $\K_n$
by finite rank operators,
and the missing factor $\#(\varphi)$ will appear along the way.

Let $T_N$ be the projection
\begin{equation}\label{eii23}
  ( T_N f) (z) = \sum_{|j|\le N } f_jz^j, \qquad N\ge 1.
\end{equation}
Note that $T_N$ is a trace class operator since it has finite
rank. Clearly $T_N, T_N^*\to 1$ strongly, and hence by Lemma
\ref{lem0} \ref{lem2},
\begin{equation}
   \det(1 -s\K_n)
= \lim_{N\to\infty} \det(1+s (P_0 -  M_\varphi P_n M_{\varphi^{-1}}) T_N).
\end{equation}

Now since $P_kT_N$ is trace class, proceeding as above in
\eqref{e10}-\eqref{e13}, we have for $N\ge |n|$,
\begin{equation}
\det(1+s (P_0 -  M_\varphi P_n M_{\varphi^{-1}}) T_N) =(1+s)^{n}
\det((1+s P_nT_N) (1-sP_n
M_{\varphi^{-1}}T_N(1-\frac{s}{1+s}P_0T_N)M_\varphi P_n )).
\end{equation}
Thus we have
\begin{equation}\label{e17}
  \det(1 -s\K_n)
= (1+s)^{n} \lim_{N\to\infty} \det(1+X_N+Y_N),
\end{equation}
where for $N\ge n$,
\begin{eqnarray}
  X_N &=& s \frac{1+sT_N}{1+s} P_n (T_N- M_{\varphi^{-1}} T_N M_\varphi ) P_n,\\
Y_N &=& -s^2 \frac{1+sT_N}{1+s} P_n M_{\varphi^{-1}} T_N (1-P_0)M_\varphi P_n.
\end{eqnarray}
We observe that
\begin{enumerate}
\item
$X_N$ and $Y_N$ are trace class.
\item $X_N^*\to 0$ strongly as $N\to\infty$.
\item $Y_N \to -s^2\S_n$ in trace norm as $N\to\infty$.
\item $(1+X_N)^{-1}$ and $(1+Y_N)^{-1}$
are uniformly bounded in operator norm as $N\to\infty$
when $s$ is small enough.
\end{enumerate}
The third property follows
using Lemma \ref{lem0} \ref{lem2}
as $(1-P_0)M_\varphi P_n$ is trace class
and $T_N\to 1$ strongly.

For a moment, we assume that $s$ is small so that
(iv) is satisfied.
Now we rewrite the right hand side of \eqref{e17} as
\begin{equation}
  \det(1+X_N+Y_N) = \det(1+Y_N)\det(1+X_N)
\det(1-\frac1{1+X_N}\frac{1}{1+Y_N}Y_NX_N).
\end{equation}
From the properties (i), (ii), (iv) above, we have
$\frac1{1+X_N}\frac{1}{1+Y_N}Y_NX_N \to 0$ in trace norm.
Using the property (iii), we now have
\begin{equation}\label{e21}
  \det(1-s\K_n) = (1+s)^{n} \det(1-s^2\S_n) \lim_{N\to\infty}
\det(1+X_N).
\end{equation}
Rewrite $X_N$ as
\begin{equation}\label{ec33}
\begin{split}
  X_N &= s \frac{1+sT_N}{1+s} P_n ( M_{\varphi^{-1}} (1-T_N) M_\varphi -(1-T_N)) P_n\\
 &= s \frac{1+sT_N}{1+s} P_n ( M_{\varphi^{-1}} P_{N+1} M_\varphi - P_{N+1}) P_n
+ s \frac{1+sT_N}{1+s} P_n ( M_{\varphi^{-1}} (1-P_{-N})M_\varphi - (1-P_{-N})) P_n\\
&=: Z_N+W_N.
\end{split}
\end{equation}
Then as $N\ge |n|$,
\begin{equation}
  W_N= s \frac{1+sT_N}{1+s} P_n M_{\varphi^{-1}} (1-P_{-N})M_\varphi P_n
= s \frac{1+sT_N}{1+s} P_n M_{\varphi^{-1}} (1-P_{-N})
(1-P_0)M_\varphi P_n,
\end{equation}
and hence $W_N\to 0$ in trace norm
as $(1-P_{-N})\to 0$ strongly and $(1-P_0)M_\varphi P_n$ is in trace class.
Also $Z_N\to 0$ strongly, and $(1+Z_N)^{-1}$, $(1+W_N)^{-1}$ are
uniformly bounded as $N\to\infty$ for $s$ small enough.
Thus by similar arguments leading to \eqref{e21},
we have
\begin{equation}\label{e24}
  \lim_{N\to\infty} \det(1+X_N)= \lim_{N\to\infty} \det(1+Z_N).
\end{equation}

Now by \eqref{e9} and \eqref{3}, we note that
\begin{equation}
  M_{\varphi^{-1}} P_{N+1} M_\varphi - P_{N+1} = \tK_N+Q_{N+1}
= A_{N+1}\ttK A_{N+1}^{-1},
\end{equation}
where $\tK_N$ is $\KN1$ with $\varphi$ replaced by $\varphi^{-1}$,
and $A_{N+1}$ is the operator of multiplication by $z^{N+1}$.
Thus we have
\begin{equation}
\begin{split}
  \det(1+Z_N) &=
\det( 1+ s \frac{1+sT_N}{1+s} P_n A_{N+1}\ttK A_{N+1}^{-1} P_n )\\
&= \det (1+ s \ttK A_{N+1}^{-1} P_n \frac{1+sT_N}{1+s} P_n A_{N+1}),
\end{split}
\end{equation}
by Lemma \ref{lem0} (iii).
Since $P_kA_{N+1}= A_{N+1}P_{k-N-1}$,
and $T_N= P_{-N}-P_{N+1}$, we have
\begin{equation}
  A_{N+1}^{-1} P_n \frac{1+sT_N}{1+s} P_n A_{N+1}
= \frac1{1+s} ((1+s)P_{n-N-1}-sP_0) \to \frac{1+s-sP_0}{1+s}=:D(s)
\end{equation}
strongly.
Since $\ttK$ is a trace class, we obtain
\begin{equation}
  \lim_{N\to\infty} \det(1+Z_N) = \det(1+s\ttK D(s)).
\end{equation}
Therefore from \eqref{e21} and \eqref{e24},
\begin{equation}\label{e29}
  \det(1-s\K_n)= (1+s)^{n}\det(1-s^2 S_n) \det(1+s\ttK D(s)).
\end{equation}

Since $\det(1+s\ttK D(s))$ does not depend on $n$, we obtain the value
of this determinant by letting $n\to\infty$ in both sides of \eqref{e29}.
But by Lemma \ref{bdj2} below,
for small $s$, $(1+s)^{-n} \det(1-s\K_n)\to (1+s)^{\#(\varphi)}$
as $n\to\infty$.
On the other hand, from \eqref{ep7}, $\det(1-s^2 S_n)$
converges to $1$.
Therefore $\det(1+s\ttK D(s))=(1+s)^{\#(\varphi)}$, and we obtain,
for small $s$,
\begin{equation}
  \det(1-s\K_n)= (1+s)^{n+\#(\varphi)}\det(1-s^2 S_n),
\end{equation}
as desired.
The result for all $s$ now follows by analytic continuation.
\end{proof}

Observe from \eqref{3}, \eqref{e9} that $\| \K_n\| \le 1+\|\varphi\|_\infty
\|\varphi^{-1}\|_\infty =: s_0$.
(A different estimate (see Appendix \cite{BDJ2}) shows that
$\|K_n\|\le \max(\|\varphi\|_\infty, \|\varphi^{-1}\|_\infty)$.)
Then we have the following result.

\begin{lem}\label{bdj2}
  For a function $\varphi$ which is
analytic and non-zero in a neighborhood of the unit circle $\{|z|=1\}$
in the complex plane, and has winding number equal to $\#(\varphi)$,
we have for $|s|<s_0^{-1}$,
\begin{equation}\label{a1}
  \lim_{n\to\infty} (1+s)^{-n}\det(1-s\K_n) = (1+s)^{\#(\varphi)},
\end{equation}
\begin{equation}\label{a2}
 \lim_{n\to -\infty} (1-s)^{n}\det(1-s\K_n) = (1-s)^{-\#(\varphi)},
\end{equation}
for some $s_0>0$.
\end{lem}
\begin{proof}
  In Lemma 5 of \cite{BDJ2}, we obtained the result \eqref{a1} for
$\varphi(z)=e^{\sqrt{\lambda}(z-z^{-1})}$ and $s=\sqrt{t}, 0\le t\le 1$.
For general analytic $\varphi$ and $s$ small, the proof remains the same
until equation (50).
The second component in the asymptotics of $F'$ is now
$-\frac{\varphi'}{1+s'}$, and hence
(51) is changed to
\begin{equation}
\begin{split}
  \log\det(1-s\K_n)
&= - \int_0^{s} \frac{ds'}{2\pi i(1+s')}
\int_{|z|=1} \biggl[-nz^{-1}- \frac{\varphi'(z)}{\varphi(z)}\biggr] dz
+ O(e^{-cn})\\
& = (n+\#(\varphi)) \log(1+s) + O(e^{-cn}).
\end{split}
\end{equation}
The calculation for \eqref{a2} is similar and we skip the details.
\end{proof}

\begin{rems}
  Lemma \ref{bdj2} does not require $\varphi$ to be analytic.
However, in this case,
the proof is particularly simple,
and can be quoted directly from \cite{BDJ2} as above.
\end{rems}
\begin{rems}
  The fact that $\det(1+s\ttK D(s))= (1+s)^{\#(\varphi)}$
is rather remarkable.
It is an instructive exercise to check this identity directly
when $\varphi$ is simple, say $\varphi=z^k$
or $\varphi=(1+az)(1+bz^{-1})$, $|a|, |b| <1$.
\end{rems}
\begin{rems}
 By \eqref{ethm2}, we see that if $n+\#(\varphi)>0$,
$\det(1-s\K_n)$ has a root at $s=-1$ of order at least
$n+\#(\varphi)$. In particular, $\K_n$ has eigenvalue $-1$.
Moreover, if $\K_n$ is self-adjoint (which is true by \eqref{e9}
whenever $|\varphi|=1$, e.g.,
$\varphi=e^{\sqrt{\lambda}(z-z^{-1})}$ as in \cite{BDJ, BDJ2}),
then $\K_n$ has an eigenspace of dimension at least
$n+\#(\varphi)$ corresponding to the eigenvalue $-1$. It is also
clear from \eqref{ethm2} that if $s\neq-1$ is a root of
$\det(1-s\K_n)$, then so is $-s$. On the other hand, if
$n+\#(\varphi)<0$, then clearly $\det(1-s^2\S_n)$ has a root at
$s=\pm 1$, etc. In the self-adjoint case, when $|\varphi|=1$, we
see from \eqref{e9} that $\S_n$ is positive definite with norm
$\le 1$. We will use this fact in Section \ref{sec-rhp}.
\end{rems}
\begin{rems}
Define the operator $A$ acting on $\ell(\Z)$ by
\begin{equation}
   A=M_{\varphi^{-1}} P_0 M_\varphi.
\end{equation}
Since $\det(1-s\K_n)=\det(1-sM_{\varphi^{-1}}\K_n M_\varphi)$,
using \eqref{e9}, the above theorem can be rephrased as
\begin{eqnarray}
  \det(1-s(P_n-A))&=& (1+s)^{n+\#(\varphi)} \det(1-s^2P_n(1-A)P_n), \\
  &=& (1-s)^{-n-\#(\varphi)} \det(1-s^2(1-P_n)A(1-P_n)).
\end{eqnarray}
These are the identities (8.55), (8.56) in \cite{Rains:corr}
for a certain subclass of $\varphi$'s with zero winding, $\#(\varphi)=0$.
\end{rems}

The following Corollary will be used in the analysis of (S3)
in Section \ref{sec-color} below.

\begin{cor}\label{cor-id}
  Let $\varphi(z)$ be as in Theorem \ref{mainthm}.
Define $\Km_n$, $\Ssm$ and $\Rrm$ to be the operators
analogous to $\K_n$, $\Ssm$ and $\Rrm$ with the matrix elements
given by
\begin{equation}
   \Km_n(z,w)
= \frac{1-z^n\varphi(z^m)w^{-n}\varphi(w^m)^{-1}}{2\pi i(z-w)},
\end{equation}
and
\begin{eqnarray}
   \Ssm(i,j)&=& \sum_{k\ge 1} (\varphi^{-1})_{(i+k)/m}
\varphi_{(-j-k)/m}, \\
   \Rrm(i,j)&=& \sum_{k\le 0} (\varphi^{-1})_{(i+k)/m}
\varphi_{(-j-k)/m},
\end{eqnarray}
where $\varphi_a=(\varphi^{-1})_a=0$ if $a\neq \Z$.
Set $\Ssmn=\chi_{[n,\infty)}\Ssm$
and $\Rrmn=\chi_{(-\infty,n-1]}\Rrm$.
Then we have
\begin{eqnarray}
  \det (1-s\Km_n) &=& (1+s)^{n+\#(\varphi)}\det (1- s^2\Ssmn),
\qquad s\neq -1,\\ &=& (1-s)^{-n-\#(\varphi)}\det (1- s^2\Rrmn),
\quad s\neq 1.
\end{eqnarray}
\end{cor}

\begin{rem}
  Observe that $\Ssmn$ has the block structure
\begin{equation}
  \begin{pmatrix}
\Ssmn(mi,mj) & \cdots & \Ssmn(mi, mj+m-1) \\
\cdots & \cdots & \cdots \\
\Ssmn(mi+m-1,mj) & \cdots & \Ssmn(mi+m-1, mj+m-1)
\end{pmatrix}
= \begin{pmatrix}
\S_n(i,j) & 0 & \cdots & 0 \\
0& \S_n(i,j) & \cdots & 0 \\
& \cdots  & &\\
0& 0& \cdots & \S_n(i,j)
\end{pmatrix}.
\end{equation}
\end{rem}

For the multi-interval case, we can generalize the argument
in Theorem \ref{mainthm} to obtain the following result.

\begin{thm}\label{thmmulti}
Let $0=n_0 \le n_1\le n_2\le \cdots \le n_k\le n_{k+1}=\infty$
be integers,
and let $s_1,\cdots, s_k$ be complex numbers satisfying
$s_k\neq -1$ and $s_k-s_j\neq -1$.
Also set $s_0=0$.
We have
\begin{equation}\label{ethmmulti}
\begin{split}
  &\det(1-\sum_{j=1}^k(s_j-s_{j-1})\Kk_{n_j})\\
&\quad
= (1+s_k)^{\#(\varphi)}
\biggl[\prod_{j=0}^{k-1}(1+s_k-s_j)^{n_{j+1}-n_j}\biggr]
\det\biggl(1-\biggl(\sum_{j=1}^{k}\frac{s_ks_j}{1+s_k-s_j}\chi_{[n_j,n_{j+1})}
\biggr)\Ss\biggr)
\end{split}
\end{equation}
where $\#(\varphi)$ is again the winding number of $\varphi$.
\end{thm}

\begin{rem}\label{rem-Bott2}
As noted by B\"ottcher (\cite{Bott2} ; ef. Remark \ref{rem-Bott1}
above), the formula \eqref{ethmmulti} remains true in the $N\times
N$ matrix case, provided we replace $n_{j+1}-n_j$ by
$N(n_{j+1}-n_j)$, $0\le j\le k-1$, and $\#(\varphi)$ by
$\#(\det\varphi)$. Again the proof in the scalar case extends to
the matrix case, and we provide no details.
\end{rem}

\begin{proof}
  The formal procedure (without considering the winding number) is as follows.
For $j=0,\cdots, k-1$, let $R_j$ be the projection operator
on $\{n_j,\cdots, n_{j+1}-1\}$,
and let $R_{k}$ be the projection operator on $\{n_k,n_k+1\cdots\}$.
Since we have from \eqref{e9}
\begin{equation}
  \Kk_{n_j}=- \sum_{l=0}^k R_l + M_\varphi \biggl( \sum_{l=j}^k
R_l \biggr) M_{\varphi^{-1}}, \qquad j=1,\cdots, k,
\end{equation}
the determinant on the left-hand-side in \eqref{ethmmulti},
denoted by $(*)$, is equal to
\begin{equation}
  (*)=\det\biggl(1+s_k \sum_{j=0}^kR_j - M_\varphi \biggl( \sum_{j=1}^k
s_j R_j \biggr) M_{\varphi^{-1}}\biggr).
\end{equation}
First we pull out the term $1+s_k \sum_{j=0}^kR_j$, then
use Lemma \ref{lem0} (iii) to obtain
\begin{equation}\label{eii54}
\begin{split}
  (*)= &\det\biggl(1+s_k \sum_{j=0}^kR_j\biggr)
\det\biggl(1-\frac1{1+s_k \sum_{j=0}^kR_j}M_\varphi \biggl( \sum_{j=1}^k
s_j R_j \biggr) M_{\varphi^{-1}}\biggr)\\
=&\det\biggl(1+s_k \sum_{j=0}^kR_j\biggr)
\det\biggl(1-\biggl( \sum_{j=1}^k s_j R_j \biggr) M_{\varphi^{-1}}
\frac1{1+s_k \sum_{j=0}^kR_j}M_\varphi \biggr).
\end{split}
\end{equation}
Now note that (recall $s_0=0$)
\begin{equation}\label{eii55}
\begin{split}
  \det\biggl(1+s_k \sum_{j=0}^kR_j\biggr)
&= \det\biggl(1+ \sum_{j=0}^{k-1} (s_k-s_j)R_j
+ \sum_{j=1}^{k} s_jR_j \biggr)\\
&= \det\biggl(1+ \sum_{j=0}^{k-1} (s_k-s_j)R_j \biggr)
\det\biggl(1+ \frac1{1+\sum_{j=0}^{k-1} (s_k-s_j)R_j}
\sum_{j=1}^{k} s_jR_j  \biggr)\\
&= \biggl[ \prod_{j=0}^{k-1} (1+s_k-s_j)^{n_{j+1}-n_j} \biggr]
\det\biggl(1+ \sum_{j=1}^{k} \frac{s_j}{1+s_k-s_j}R_j  \biggr).
\end{split}
\end{equation}
Using \eqref{eii55} and then multiplying two determinants,
we have
\begin{equation}\label{eii56}
\begin{split}
(*)&= \biggl[ \prod_{j=0}^{k-1} (1+s_k-s_j)^{n_{j+1}-n_j} \biggr]\\
&\qquad \times \det\biggl(1+ \sum_{j=1}^{k} \frac{s_j}{1+s_k-s_j}R_j
-(1+s_k) \biggl( \sum_{j=1}^k \frac{s_j}{1+s_k-s_j} R_j \biggr)
M_{\varphi^{-1}}
\frac1{1+s_k \sum_{j=0}^kR_j}M_\varphi \biggr).
\end{split}
\end{equation}
Finally, using
\begin{equation}
  \frac1{1+s_k \sum_{j=0}^kR_j}
=\frac{1+s_k(1-\sum_{j=0}^kR_j)}{1+s_k}
\end{equation}
in the determinant on the right-hand-side of \eqref{eii56},
we obtain
\begin{equation}
(*)= \biggl[ \prod_{j=0}^{k-1} (1+s_k-s_j)^{n_{j+1}-n_j} \biggr]
  \det\biggl( 1-\sum_{j=0}^{k}\frac{s_ks_j}{1+s_k-s_j} R_j M_{\varphi^{-1}}
\biggl(1-\sum_{j=0}^kR_j\biggr) M_\varphi \biggr),
\end{equation}
which is precisely \eqref{ethmmulti} from \eqref{e9}.

The rigorous proof is also similar to the proof of Theorem \ref{mainthm}.
Let $T_N$ be the projection on $|j|\le N$ as in \eqref{eii23}.
We take $N$ large so that $N>n_k$.
The analogue of \eqref{e17} is now
\begin{equation}
  (*)=\biggl[ \prod_{j=0}^{k-1} (1+s_k-s_j)^{n_{j+1}-n_j} \biggr]
\lim_{N\to\infty} \det(1+X_N+Y_N),
\end{equation}
where
\begin{align}
   X_N &= \biggl( \sum_{j=1}^k \frac{s_j}{1+s_k-s_j}
\frac{1+s_k-s_j(1-T_N)}{1+s_k}R_j \biggr)
(T_N-M_{\varphi^{-1}}T_N M_\varphi),\\
  Y_N &= -\biggl(\sum_{j=1}^k \frac{s_ks_j}{1+s_k-s_j}\frac{1+s_k-s_j(1-T_N)}
{1+s_k}R_j \biggr)
M_{\varphi^{-1}}T_N \bigl(1-\sum_{j=0}^k s_jR_j\bigr)M_\varphi,
\end{align}
which becomes, by the same argument leading to \eqref{e24},
\begin{equation}
  (*)=\biggl[ \prod_{j=0}^{k-1} (1+s_k-s_j)^{n_{j+1}-n_j} \biggr]
\det\biggl(1-\biggl(\sum_{j=1}^{k}\frac{s_ks_j}{1+s_k-s_j}\chi_{[n_j,n_{j+1})}
\biggr)\Ss\biggr)
\lim_{N\to\infty} \det(1+Z_N),
\end{equation}
with $Z_N$ in \eqref{ec33} where $s$ is replaced by $s_k$,
This then leads to the desired result as in the single
interval case.
%and $R_n$ is replaced by $R_k$.
%\begin{equation}
%  Z_N= \biggl( \sum_{j=1}^k \frac{s_j}{1+s_k-s_j}
%\frac{1+s_kT_N}{1+s_k}R_j \biggr)
%(M_{\varphi^{-1}}R^{(N+1)} M_\varphi-R^{(N+1)}),
%\end{equation}
%where $R^{(N+1)}$ denotes the projection on $j\ge N+1$
%which was denoted by $R_{N+1}$ in the proof of Theorem \ref{mainthm}.
\end{proof}

\section{Convergence of moments}\label{sec-mom}

%Let $Y_N$ denote the set of Young diagrams 
%$\lambda=(\lambda_1,\lambda_2,\cdots)$ 
%of size $N$.
%Here $\lambda_j$ is the number of boxes in 
%the $j$th row.
%Note that $\lambda_1\ge\lambda_2\ge\cdots$.
%We equip $Y_N$ with \emph{Plancherel measure} : 
%\begin{equation}
%  \Prob^{\Plan}_N(\lambda) = \frac{d_{\lambda}^2}{N!}, 
%\qquad \lambda \in Y_N.
%\end{equation} 
%Set 
%\begin{equation}\label{en2}
%  \xi_j:=\frac{\lambda_j-2\sqrt{N}}{N^{1/6}}.
%\end{equation}

In this section, we prove the convergence of moments for 
arbitrary (scaled) rows, $\xi_j$, of a random young diagram under the 
Plancherel measure, mentioned in the Introduction.
The tail estimates used in the proof of Theorem \ref{thmmom} 
are given in Section \ref{sec-tail} below.

%converges in distribution to the \emph{same} limiting random variables : 
%for any $k\in\N$, 
%\begin{equation}\label{eq2.5}
%  \lim_{N\to\infty}\Prob^{\Plan}_N\bigl(\xi_1\le t_1, \cdots, 
%\xi_k\le t_k\bigr)
%= F(t_1,\cdots,t_k).
%\end{equation} 
%The case when $k=1$ was proved in \cite{BDJ}, 
%and the case when $k=2$ with $t_1=\infty$ was proved in \cite{BDJ2}.
%The general case was conjectured in the above two papers, and 
%was proved in \cite{Ok, BOO, kurtj:disc}.
%
%In addition to the convergence in distribution, in \cite{BDJ, BDJ2}, 
%the convergence of all the moments are proved for the corresponding cases 
%$k=1,2$.
%In this paper, we establish the convergence of all the moments 
%in the general case.

Let $\N_0:=\N\cup\{0\}$.

\begin{thm}\label{thmmom}
For any fixed $k\in\N$, and for any $a_j\in\N_0$, $1\le j\le k$, 
we have as $N\to\infty$
\begin{equation}\label{ec3.5}
  \Exp^{\Plan}_N\bigl(\xi_1^{a_1}\cdots \xi_k^{a_k}\bigr) \to 
\Exp\bigl(x_1^{a_1}\cdots x_k^{a_k}\bigr)
\end{equation} 
where $\Exp^{\Plan}_N$ denotes the expectation with respect to the 
Plancherel measure on $Y_N$, 
and $\Exp$ denotes the expectation with respect to 
the limiting distribution function $F$ 
in \eqref{eq2.4}, \eqref{ei5}.
\end{thm}

\begin{rem}\label{rem3.2}
It will be clear from the proof below that the following 
stronger convergence result is also true : 
Let $h_j(x)$, $j=1,\cdots, k$ be continuous functions on $\R$
satisfying $|h_j(x)|\le C_1e^{c_2|x|^{3/2-\epsilon}}$
for some $\epsilon>0$. Then for any $k$, as $N\to\infty$, 
\begin{equation}\label{ewq}
  \Exp^{\Plan}_N\bigl(h_1(\xi_1)\cdots h_k(\xi_k)\bigr) \to 
\Exp\bigl(h_1(x_1)\cdots h_k(x_k)\bigr).
\end{equation} 
\end{rem}

\begin{proof}
%It is enough to prove for even integers $a_j$.
%Hence we will assume that $a_j$ are all even.
%Let $P_N(\xi_1,\cdots,\xi_k)$ be the distribution function 
%for the scaled random variables $\xi_1, \cdots, \xi_k$. 
%Since $\lambda=(\lambda_1,\cdots)$, 
%we will use the notation $\Prob^{\Plan}_N$ 
%for the distribution function for $\lambda_j$.
We have 
\begin{equation}\label{e2.3}
   \Exp^{\Plan}_N(\xi_1^{a_1}\cdots \xi_k^{a_k})
= \int_{x_1\ge\cdots\ge x_k} 
\prod_{j=1}^k x_j^{a_j} d \Prob^{\Plan}_N(\xi_1\le x_1,\cdots,\xi_k\le x_k)
\end{equation} 
since $\lambda_1\ge\lambda_2\cdots$.
Fix a number $T>2$.
We split the integral into two pieces : 
\begin{align}
\label{eq3.8}
(a) & \max_{1\le j\le k}{|x_j|}\le T \\
\label{eq3.9}
(b) & \max_{1\le j\le k}{|x_j|}> T.
\end{align}
In the first part (a), using a standard argument
and the convergence in distribution \eqref{ei5} above, 
the limit becomes as $N\to\infty$
\begin{equation}\label{e2.6}
   \int_{\substack{x_1\ge\cdots\ge x_k \\ \max{|x_j|}\le T}} 
\prod_{j=1}^k x_j^{a_j} d F(x_1,\cdots,x_k).
\end{equation} 

For the second part (b), 
the region is a union of two (not necessarily disjoint) pieces : 
\begin{align}
(i) \ \ & \max_{j}{|x_j|}=|x_1|\\
(ii) \ \ & \max_{j}{|x_j|}=|x_k|.
\end{align}
Note that since $x_1\ge\cdots\ge x_k$,
$\max_{j}{|x_j|}$ is either $|x_1|$ or $|x_k|$.
Over region (i), 
\begin{equation}\label{e2.10}
\begin{split}
  \int_{(i)} \prod_{j=1}^k |x_j|^{a_j} 
d \Prob^{\Plan}_N(\xi_1\le x_1,\cdots,\xi_k\le x_k)
&\le \int_{(i)} |x_1|^{a_1+\cdots+a_k} 
d \Prob^{\Plan}_N(\xi_1\le x_1,\cdots,\xi_k\le x_k)\\
&\le \int_{(-\infty,-T)\cup (T,\infty)} |x_1|^{a_1+\cdots+a_k} 
d \Prob^{\Plan}_N(\xi_1\le x_1)\\
&= \Exp^{\Plan}_N(|\xi_1|^{a_1+\cdots+a_k}(\chi_{\xi_1<-T}+\chi_{\xi_1>T})).
\end{split}
\end{equation} 
Similarly, 
\begin{equation}\label{e2.11}
  \int_{(ii)} \prod_{j=1}^k |x_j|^{a_j} 
d \Prob^{\Plan}_N(\xi_1\le x_1,\cdots,\xi_k\le x_k)
%\le \int_{(-\infty,-T)} \xi_k^{a_1+\cdots+a_k} d \Prob^{\Plan}_N(\xi_k)
\le \Exp^{\Plan}_N(|\xi_k|^{a_1+\cdots+a_k}(\chi_{\xi_k<-T}+\chi_{\xi_k>T})).
\end{equation}

Now from the tail estimates in Proposition \ref{prop3.1} below, 
the moment \eqref{e2.3} as $N\to\infty$ is equal to \eqref{e2.6} 
plus a term which can be made arbitrarily small if 
we take $T$ large enough.
However, from Lemma \ref{lem2.2}, for $T$ large,  
\eqref{e2.6} is arbitrarily close to 
$\Exp(x_1^{a_1}\cdots x_k^{a_k})$.
Thus we have proved the theorem.
\end{proof}

\begin{lem}\label{lem2.2}
For any $k\in\N$, and for any $a_j\in\N_0$, $1\le j\le k$, 
\begin{equation}
  \Exp(|x_1|^{a_1}\cdots |x_k|^{a_k}) <\infty, 
\end{equation} 
where $\Exp$ is the expectation with respect to 
the limiting distribution function $F$ in  
\eqref{eq2.4} and \eqref{ei5}.
\end{lem}

\begin{proof}
We need to show that 
\begin{equation}\label{eq3.18}
%  \Exp(x_1^{a_1}\cdots x_k^{a_k}) =
 \int_{x_1\ge\cdots\ge x_k} 
\prod_{j=1}^k |x_j|^{a_j} d F(x_1,\cdots,x_k) <\infty.
\end{equation} 
Fix $T>2$.
We split the integral into two parts as in \eqref{eq3.8}, \eqref{eq3.9} : 
(a) $\max_j |x_j| \le T$, and (b) $\max_j |x_j| > T$.
In (a), the integral is finite.
In (b), the argument yielding \eqref{e2.10}, \eqref{e2.11} 
implies that 
\begin{equation}\label{eq3.20}
 \int_{(b)}
\prod_{j=1}^k |x_j|^{a_j} d F(x_1,\cdots,x_k) 
\le \Exp(|x_1|^{a_1+\cdots +a_k} \chi_{x_1>T})
+ \Exp(|x_k|^{a_1+\cdots +a_k} \chi_{x_k<-T}).
\end{equation} 
(Note that the additional terms corresponding to 
$\Exp^{\Plan}_N(|\xi_1|^{a_1+\cdots+a_k}\chi_{\xi_1<T})$ 
in \eqref{e2.10} 
and $\Exp^{\Plan}_N(|\xi_k|^{a_1+\cdots+a_k}\chi_{\xi_k>T})$ 
in \eqref{e2.11}
are not necessary here as $dF$ is a smooth measure.)
We will prove the finiteness of the last two expected values 
for $a=a_1+\cdots +a_k$.

First, we prove that
$\Exp(x_1^{a} \chi_{x_1>T})<\infty$ for any $a\in \N_0$.
Note that by \eqref{ei5} and \eqref{e3.1} below, 
for $x_1>T_0$, 
\begin{equation}
 1- F(x_1)
=\lim_{N\to\infty} 1-\Prob_N^{\Plan}(\xi_1\le x_1)
= \lim_{N\to\infty} \Exp_N^{\Plan} (\chi_{\xi>x_1})
\le Ce^{-cx_1^{3/2}}
\end{equation} 
for some $C,c>0$.
In particular, we have for any $a\in \N_0$, 
\begin{equation}
  \lim_{x_1\to\infty} x_1^a (1- F(x_1)) =0.
\end{equation} 
Thus, integrating by parts, 
\begin{equation}
  \int_{T}^\infty x_1^a dF(x_1)
= T^a(1-F(T))
+ \int_T^\infty ax_1^{a-1}(1-F(x_1)) dx_1.
\end{equation} 
Using \eqref{ei5}, Fatou's lemma and \eqref{e3.1}, we have 
\begin{equation}
\begin{split}
   \Exp(x_1^a \chi_{x_1>T})
  &=\int_{T}^\infty x_1^a dF(x_1) \\
& = \lim_{N\to\infty} T^a(1-\Prob^{Plan}_N(\xi_1\le T)) 
+ \int_T^\infty \lim_{N\to\infty} 
ax_1^{a-1}(1-\Prob^{Plan}_N(\xi_1\le x_1)) dx_1\\
&\le \liminf_{N\to\infty}
\biggl[   T^a(1-\Prob^{Plan}_N(\xi_1\le T)) 
+ \int_T^\infty ax_1^{a-1}(1-\Prob^{Plan}_N(\xi_1\le x_1)) dx_1 \biggr]\\
&= \liminf_{N\to\infty} \Exp_N^{\Plan} (\xi_1^{a} \chi_{\xi_1>T})\\
&\le Ce^{-cT^{3/2}} <\infty.
\end{split}
\end{equation} 

The proof of the finiteness of the second expected value 
in \eqref{eq3.20} is similar using \eqref{e3.2}.
\end{proof}

\section{Tail estimates}\label{sec-tail}

%Recall \eqref{en2}.

For the proof of Theorem \ref{thmmom}, we need tail estimates 
for the (scaled) length $\xi_k$ of each row, which are uniform in $N$.
In this section, we obtain these tail estimates 
in Proposition \ref{prop3.1}.
These estimates follow from the tail estimates, Proposition \ref{prop3.2} 
for the Poissonized 
Plancherel measure introduced in Section \ref{sec-intro}, 
together with the de-Poissonization Lemma \ref{lem3.3}.

%As mentioned in the Introduction, %in \cite{Jo2, BDJ, BDJ2, BOO, kurtj:disc}, 
%we work with the \emph{Poissonized Plancherel measure} $\Prob^{Pois}_t$, 
%rather than working directly with $\Prob^{Plan}_N$. 
Define 
\begin{equation}\label{ewr}
  \phi^{(k)}_n(t):= 
\Prob^{Pois}_t(\lambda_k\le n)
=
\sum_{N=0}^\infty 
\frac{e^{-t^2}t^{2N}}{N!} \Prob^{\Plan}_N(\lambda_k\le n).
\end{equation} 
(In \cite{BDJ, BDJ2}, the notation $\lambda=\sqrt{t}$ is used. 
But in this paper, to avoid the confusion with 
the notation $\lambda$ for a partition, 
we use $t$.)
The following result is proved in Section \ref{sec-rhp} using the 
steepest-descent method for RHP.
%We first obtain tail estimations for $\phi^{(k)}_n(t)$.
%As in \cite{Jo2, BDJ}, a de-Poissonization lemma 
%(see Lemma \ref{lem3.3} below) will imply tail 
%estimates for the Plancherel measure from estimates for 
%the Poissonized Plancherel measure.
Note that $0\le \phi^{(k)}_n(t)\le 1$.

\begin{prop}\label{prop3.2}
 Define $x$ by 
\begin{equation}\label{e3.4}
  \frac{2t}{n}=1-\frac{x}{2^{1/3}n^{2/3}}.
\end{equation} 
Let $k\in \N$.
There are constants $C,c>0$ and $0<\delta_0<1$ such that 
for large $t$ and $n$, and for any fixed $0<\delta<\delta_0$, 
the following hold true : 
for $x\ge 0$, 
\begin{eqnarray}
\label{en3}
  &0\le 1- \phi^{(k)}_n(t) \le C e^{-cn}, \quad 
&0\le \frac{2t}{n}\le 1-\delta, \\
\label{e3.5}
  &0\le 1- \phi^{(k)}_n(t) \le C e^{-cx^{3/2}}, \quad 
&1-\delta<\frac{2t}{n}\le 1, 
\end{eqnarray} 
and for $x<0$, 
\begin{eqnarray}
\label{e3.6}
  &0\le  \phi^{(k)}_n(t) \le Ce^{-c |x|^{3/2}}, \quad
& 1< \frac{2t}{n}< 1+\delta, \\
\label{en4}
  &0\le  \phi^{(k)}_n(t) \le Ce^{-c t}, \quad
& 1+\delta\le \frac{2t}{n}.
\end{eqnarray} 
\end{prop}

%The proof of this Proposition will be given in Section \ref{sec-rhp} 
%below.
%To prove Proposition \ref{prop3.1}, 

We also need the following de-Poissonization lemma : 

\begin{lem}\label{lem3.3}
  There exists $C>0$ such that for all sufficiently large $N$, 
\begin{equation}
  \Prob^{\Plan}_N(\lambda_k\le n) \le 
C\phi^{(k)}_n\bigl((N-\sqrt{N})^{1/2}\bigr),
\qquad 
1- \Prob^{\Plan}_N(\lambda_k\le n) \le 
C\bigl(1-\phi^{(k)}_n\bigl((N+\sqrt{N})^{1/2}\bigr)\bigr)
\end{equation} 
for all $n\in\Z$.
\end{lem}

\begin{proof}
  This is similar to Lemma 8.3 in \cite{BDJ} 
(again note that $\lambda$ in \cite{BDJ} satisfies $\lambda=\sqrt{t}$.) 
Indeed, the proof of Lemma 8.3 in \cite{BDJ} only requires 
the fact that $0\le q_{n,N+1}\le q_{n,N}\le 1$.
In our case, $q_{n,N}=\Prob^{\Plan}_N(\lambda_k\le n)$, 
which is clearly between $0$ and $1$.
The monotonicity can be found in \cite{kurtj:disc} Lemma 3.8.
\end{proof}

Now Proposition \ref{prop3.2} and Lemma \ref{lem3.3} imply the following 
uniform tail estimates.

\begin{prop}\label{prop3.1}
  Fix $k\in \N$ and $a\in \N_0$.
For a given $T\ge 2$, 
there are constants $C, c>0$ and $N_0>0$ such that 
for $N\ge N_0$, 
\begin{equation}\label{e3.1}
  \Exp^{\Plan}_N(\xi_k^{a}\chi_{\xi_k>T})\le Ce^{-cT^{3/2}}
+ Ce^{-cN^{1/2}}
\end{equation} 
and 
\begin{equation}\label{e3.2} 
\Exp^{\Plan}_N(|\xi_k|^{a}\chi_{\xi_k<-T}) 
\le Ce^{-cT^{3/2}} + Ce^{-cN^{1/2}}.
\end{equation} 
\end{prop}

\begin{proof}  %[Proof of Proposition \ref{prop3.1}]

{\bf (a) Bound \eqref{e3.1} :} 
%We start with the bound \eqref{e3.1}.
%Since $\xi_1\ge\xi_k$, $\xi_k>T$ implies $\xi_1>T$. 
%Hence 
%\begin{equation}
%  \Exp^{\Plan}_N(\xi_k^{a}\chi_{\xi_k>T})
%\le \Exp^{\Plan}_N(\xi_1^{a}\chi_{\xi_1>T}).
%\end{equation} 
%Thus it is enough to prove \eqref{e3.1} for $k=1$.
Without any loss we can assume $a>0$.
Note that since $0\le\lambda_k\le N$, 
\begin{equation}\label{e3.9}
  -2N^{1/3}\le \xi_k\le \frac{N-2\sqrt{N}}{N^{1/6}}< N^{5/6}.
\end{equation} 
If $T\ge N^{5/6}$, then the expected value in \eqref{e3.1} 
is zero, and the bound is trivial.
Thus we assume that $T< N^{5/6}$. 
Integrating by parts and using Lemma \ref{lem3.3},
\begin{equation}\label{e3.10}
\begin{split}
   \Exp^{\Plan}_N(\xi_k^{a}\chi_{\xi_k>T})
&= \int_{(T,N^{5/6})} s^a d\Prob^{\Plan}_N(\xi_k\le s) \\
%\end{equation} 
%Integrating by parts, and using Lemma \ref{lem3.3}
%and Proposition \ref{prop3.2}, 
%for some constants $C, c>0$,
%\begin{equation}\label{e3.10}
%\begin{split}
%   \Exp^{\Plan}_N(\xi_k^{a}\chi_{\xi_k>T})
&= T^a\bigl(1-\Prob^{\Plan}_N(\xi_k\le T)\bigr)
+ \int_{(T,N^{5/6})} as^{a-1}\bigl(1-\Prob^{\Plan}_N(\xi_k\le s)\bigr) ds\\
&\le CT^a(1-\phi_{n(T)}^{(1)}((N+\sqrt{N})^{1/2}))
+ C\int_{(T,N^{5/6})} as^{a-1}(1-\phi_{n(s)}^{(1)}((N+\sqrt{N})^{1/2})) 
ds,
%&\le CT^ae^{-cx(T,N)^{3/2}} 
%+ \int_{(T,N^{5/6})} Cas^{a-1} e^{-cx(s,N)^{3/2}} ds,
\end{split}
\end{equation} 
for large $N$, 
where $n(s)=2\sqrt{N}+sN^{1/6}$.
Note that since $T\ge 2$, 
$\frac{2(N+\sqrt{N})^{1/2}}{n(T)}\le 1$.
We distinguish two cases : 
\begin{eqnarray}
  &(i) \quad 0\le \frac{2(N+\sqrt{N})^{1/2}}{n(T)}< 1-\delta \\
&(ii) \quad 1-\delta \le \frac{2(N+\sqrt{N})^{1/2}}{n(T)}\le 1
\end{eqnarray}
where $0<\delta<1$ is a fixed constant satisfying $\delta <\delta_0$,
where $\delta_0$  
appears in Proposition \ref{prop3.2}.

{\bf Case (i) :}
For all $s\ge T$, $0\le \frac{2(N+\sqrt{N})^{1/2}}{n(s)}\le 
\frac{2(N+\sqrt{N})^{1/2}}{n(T)} < 1-\delta$.
Note that for $T\le s$, 
\begin{equation}
  n(s)\ge n(T)\ge \frac{2(N+\sqrt{N})^{1/2}}{1-\delta}
\ge \frac{2\sqrt{N}}{1-\delta}.
\end{equation}
Hence from the estimate \eqref{en3}, we have 
\begin{equation}
  1-\phi^{(1)}_{n(s)} ((N+\sqrt{N})^{1/2})
\le Ce^{-cn(s)} \le Ce^{-cN^{1/2}}
\end{equation}
for $T\le s<N^{5/6}$ with a new constant $c$.
Therefore, from \eqref{e3.10}, we obtain
\begin{equation}
  \Exp_N^{\Plan} (\xi_k^a \chi_{\xi_k>T}) \le Ce^{-cN^{1/2}}.
\end{equation}

{\bf Case (ii) : }
There is $s_0>T$ such that $\frac{2(N+\sqrt{N})^{1/2}}{n(s_0)} 
=1-\delta$.
We write \eqref{e3.10} as 
\begin{equation}
\begin{split}
  \Exp_N^{\Plan} (\xi_k^a \chi_{\xi_k>T}) 
&\le CT^a(1-\phi^{(1)}_{n(T)} ((N+\sqrt{N})^{1/2})) 
+ C\int_{(T,s_0)} as^{a-1} 
(1-\phi^{(1)}_{n(s)} ((N+\sqrt{N})^{1/2})) ds  \\
&\quad + C\int_{[s_0,N^{5/6})} as^{a-1} 
(1-\phi^{(1)}_{n(s)} ((N+\sqrt{N})^{1/2})) ds 
\qquad \text{(if $s_0\ge N^{5/6}$, the third integral is zero)} \\
&\le CT^ae^{-cx(T)^{3/2}} + C\int_{(T,s_0)} as^{a-1} e^{-cx(s)^{3/2}}ds
+ C\int_{[s_0,N^{5/6})} as^{a-1} e^{-cn(s)}ds
\end{split}
\end{equation}
using \eqref{en3}, \eqref{e3.5}, 
where $x(s)$ is defined by the formula \eqref{e3.4} 
with $t=(N+\sqrt{N})^{1/2}$ and $n=2\sqrt{N}+sN^{1/6}$.
As in Case (i), for $s\ge s_0$, 
\begin{equation}
  n(s)\ge n(s_0) = \frac{2(N+\sqrt{N})^{1/2}}{1-\delta}
\ge \frac{2\sqrt{N}}{1-\delta}, 
\end{equation}
and hence, the last integral is less than $Ce^{-cN^{1/2}}$.
For the other terms, since 
\begin{equation}
  2\sqrt{N}-2(N+\sqrt{N})^{1/2} 
=\frac{-4\sqrt{N}}{2\sqrt{N}+2(N+\sqrt{N})^{1/2}} \ge -1
\ge -N^{1/6}, 
\end{equation}
we have for $T\le s<s_0$, 
\begin{equation}
  x(s) = \frac{n-2t}{(n/2)^{1/3}}
=\frac{2\sqrt{N}+sN^{1/6}-2(N+\sqrt{N})^{1/2}}
{(\sqrt{N}+\frac{s}2N^{1/6})^{1/3}}
\ge \frac{s-1}
{(1+\frac{s}{2N^{1/3}})^{1/3}}
\ge \frac{\frac12}
{(1+\frac{s}{2N^{1/3}})^{1/3}}
\end{equation}
as $s\ge T\ge 2$.
Noting that 
$s_0= \frac{2(N+\sqrt{N})^{1/2}-2(1-\delta)\sqrt{N}}{(1-\delta)N^{1/6}} 
\le c_0N^{1/3}$ for some constant $c_0$, 
we have 
\begin{equation}
   x(s)\ge \frac{\frac12}
{(1+\frac{s}{2N^{1/3}})^{1/3}} \ge cs
\end{equation}
for $s\le s_0$ 
with some constant $c>0$.
Hence 
\begin{equation}
\begin{split}
  \Exp^{\Plan}_N (\xi_k^a \chi_{\xi_k>T}) 
& \le CT^ae^{-cT^{3/2}} + C\int_T^\infty s^{a-1}e^{-cs^{3/2}} ds
+ Ce^{-cN^{1/2}} \\
&\le Ce^{-cT^{3/2}} + Ce^{-cN^{1/2}}.
\end{split}
\end{equation}

{\bf (b) Bound \eqref{e3.2} :}
%Now we prove the bound \eqref{e3.2}.
Recalling \eqref{e3.9}, 
if $T\ge 2N^{1/3}$, the expected value in \eqref{e3.2} 
is zero and the bound is trivial.
Thus we assume that $T\le 2N^{1/3}$.
%Using \eqref{e3.9}, Lemma \ref{lem3.3} and \ref{prop3.2}, 
Integrating by parts and using Lemma \ref{lem3.3}, 
we have for some constants $C,c>0$,
\begin{equation}\label{e3.19}
\begin{split}
   \Exp^{\Plan}_N(|\xi_k^{a}|\chi_{\xi_k<-T})
&= \int_{[-2N^{1/3},-T)} (-s)^a d \Prob^{\Plan}_N(\xi_k\le s)\\
&= |-T|^a\Prob^{\Plan}_N(\xi_k<- T)
+ \int_{[-2N^{1/3},-T)} a(-s)^{a-1}\Prob^{\Plan}_N(\xi_k\le s) ds\\
&\le CT^a\phi^{(k)}_{n(-T)}((N-\sqrt{N})^{1/2}
+ C\int_{[-2N^{1/3},-T)}a(-s)^{a-1}\phi^{(k)}_{n(s)}((N-\sqrt{N})^{1/2}
ds
%&\le CT^ae^{-c|x(-T,N)|^{3/2}}
%+ \int_{[-2N^{1/3},-T)} Ca|s|^{a-1} e^{-c|x(s,N)|^{3/2}} ds,
\end{split}
\end{equation} 
for large $N$, where $n(s)=2\sqrt{N}+sN^{1/6}$ as before.
Given $T$, we take $N_0>\frac{1+\sqrt{1+T^3}}{2T}$ 
so that for $N\ge N_0$, $\frac{2(N-\sqrt{N})^{1/2}}{n(-T)}>1$.
We distinguish two cases : 
\begin{eqnarray}
  &(i) &1+\delta \le \frac{2(N-\sqrt{N})^{1/2}}{n(-T)} \\
  &(ii) &1< \frac{2(N-\sqrt{N})^{1/2}}{n(-T)} <1+\delta
\end{eqnarray}
where $0<\delta<1$ is a fixed constant as above.

{\bf Case (i) :}
For all $-2N^{1/3}-2N^{1/3}<s<-T$, $\frac{2(N-\sqrt{N})^{1/2}}{n(s)} 
\ge \frac{2(N-\sqrt{N})^{1/2}}{n(-T)} \ge 1+\delta$.
From the estimate \eqref{en4}, using $T\le 2N^{1/3}$, 
\begin{equation}
  \Exp^{\Plan}_N(|\xi^a_k|\chi_{\xi_k<-T}) \le Ce^{-cN^{1/2}}.
\end{equation}

{\bf Case (ii) :}
There is $s_0>T$ such that 
$\frac{2(N-\sqrt{N})^{1/2}}{n(-s_0)}=1+\delta$.
We write \eqref{e3.19} as 
\begin{equation}
\begin{split} 
   \Exp^{\Plan}_N(|\xi^a_k|\chi_{\xi_k<-T}) 
&\le CT^a\phi^{(k)}_{n(-T)}((N-\sqrt{N})^{1/2})
+ C \int_{(-s_0,-T)} a(-s)^{a-1} 
\phi^{(k)}_{n(-T)}((N-\sqrt{N})^{1/2})ds \\
&\quad + C\int_{[-2N^{1/3},-s_0]} a(-s)^{a-1} 
\phi^{(k)}_{n(-T)}((N-\sqrt{N})^{1/2}) ds \\
&\le CT^ae^{-c|y(-T)|^{3/2}} 
+ C\int_{(-s_0,-T)} a(-s)^{a-1} e^{-c|y(s)|^{3/2}} ds \\
&\quad + C\int_{[-2N^{1/3},-s_0]} a(-s)^{a-1} 
e^{-c(N-\sqrt{N})^{1/2}}ds
\end{split} 
\end{equation} 
using \eqref{e3.6}, \eqref{en4}, 
where $y(s)$ is defined by $x$ in the formula \eqref{e3.4}
with $t=(N-\sqrt{N})^{1/2}$ and $n=2\sqrt{N}+sN^{1/6}$.
The last integral is less than $Ce^{-cN^{1/2}}$ with 
new constants $C,c$.
For the other terms, note that 
since $2\sqrt{N}-2(N-\sqrt{N})^{1/2}\le 2\le N^{1/6}$ and $s\le -T\le -2$, 
we have for $-2N^{1/3}\le s\le -T$, 
\begin{equation}
  y(s)= \frac{n-2t}{(n/2)^{1/3}}
= \frac{s+ [2\sqrt{N}-2(N-\sqrt{N})^{1/2}]N^{-1/6}}
{(1+s/(2N^{1/3}))^{1/3}}.
\le s+1 \le \frac{s}2. 
\end{equation}
Thus, we obtain 
\begin{equation}
\begin{split} 
  \Exp^{\Plan}_N(|\xi^a_k|\chi_{\xi_k<-T})
&\le CT^ae^{-cT^{3/2}}
+ C\int_{(-s_0,-T)} a|s|^{a-1} e^{-c|s|^{3/2}} ds
+ Ce^{-cN^{1/2}} \\
&\le Ce^{-cT^{3/2}} + Ce^{-cN^{1/2}}.
\end{split} 
\end{equation} 
\end{proof}
 
\begin{rem}
The results \eqref{en3} - \eqref{en4} 
for $k=1$ were given in \cite{BDJ}.
Indeed, in \cite{BDJ} stronger bounds than \eqref{e3.6}, \eqref{en4} 
were obtained (Lemma 7.1 (iv), (v) in \cite{BDJ}) :
\begin{eqnarray} 
\label{q1}
  &0\le  \phi^{(1)}_n(t) \le Ce^{-c |x|^{3}}, \quad 
& 1<\frac{2t}{n} <1+\delta, \\
\label{q2}
  &0\le  \phi^{(1)}_n(t) \le Ce^{-c t^2}, \quad 
&  
1+\delta \le \frac{2t}{n}
\end{eqnarray} 
(note $\lambda=\sqrt{t}$ in \cite{BDJ}.)
From this, as in the Proof of Proposition \ref{prop3.1}. 
we have 
\begin{equation}\label{e3.3}
  \Exp^{\Plan}_N(|\xi_1^{a}|\chi_{\xi_1<-T})
\le Ce^{-cT^3}+ Ce^{-cN}.
\end{equation}
In this paper, we only obtain the above weaker bounds 
\eqref{e3.6} and \eqref{en4}, but
they are enough for our purpose in proving the convergence of moments.
However, we believe that the same bound \eqref{e3.3} holds true for
general $k$.
In the next section, we indicate why we only obtain these weaker bounds
(see the Remark before Lemma \ref{lem5.1} below).
\end{rem}

\section{Riemann-Hilbert problem}\label{sec-rhp}

In this section, we will prove Proposition \ref{prop3.2}.

For \eqref{en3} and \eqref{e3.5}, note that 
$\phi^{(k)}_n(t)\ge \phi^{(1)}_n(t)$ for all $k\ge 1$
as $\lambda_1\ge\lambda_2\ge\cdots$.
But for $k=1$ the estimates \eqref{en3}, \eqref{e3.5} 
were proved in \cite{BDJ} (Lemma 7.1 (i), (ii)), 
and so we have the same bounds for all $k\ge 1$.
On the other hand, since $\Prob^{\Plan}_N(\lambda\le n)\le 1$, 
we always have $0\le \phi^{(k)}_n(t)\le 1$.

%We will focus on \eqref{e3.6} in t
The rest of this section is devoted to proving \eqref{e3.6} 
and \eqref{en4}.
We start from the formulae (see \cite{BOO, kurtj:disc, Ok99, Rains:corr}) 
that (recall $\phi^{(k)}_n(t)=\Prob^{Pois}_t(\lambda_k\le n)$)
\begin{eqnarray}
   \phi^{(1)}_n(t)&=& \det(1-\S_n),\\
\label{ez1}   \phi^{(k+1)}_{n+k}(t)&=& \phi^{(k)}_{n+k-1}(t) + 
\biggl(-\frac{d}{dr}\biggr)^k\bigg|_{r=1} \det(1-r\S_n), 
\qquad k\ge 1,
\end{eqnarray} 
with 
\begin{equation}\label{e5.3}
   \varphi(z)=e^{t(z-z^{-1})}.
\end{equation} 
This follows from, for example, Theorem 3.1 of \cite{Rains:corr} 
with $p_+=p_-=(t,0,0)$ (see also \cite{Rains:mean}) which states that 
for any finite subset $A$ of $\Z$, 
\begin{equation}
  \Prob^{Pois}_t( A\subset \{ \lambda_j-j\} ) = \det ( S(i,j))_{i,j\in A}, 
\qquad S(i,j) = \sum_{k=1}^\infty (\varphi^{-1})_{i+k}\varphi_{-j-k},
\end{equation} 
with $\varphi$ given by \eqref{e5.3},
where $\Prob^{Pois}_t$ 
denotes the Plancherel measure for $S_N$ with $N$ being the 
Poisson variable with mean $t^2$ (Poissonized Plancherel measure).
Recall from \eqref{ewr} 
that $\phi^{(k)}_{n}(t)=\Prob^{Pois}_t(\lambda_k\le n)$.
In \cite{BDJ}, the authors obtained the estimates (stronger than) 
\eqref{e3.6} and \eqref{en4} in the case $k=1$ (see the Remark at the end 
of Section \ref{sec-tail}).
% $\det(1-\S_n)\le Ce^{-c|x|^{3/2}}$, $x\le 0$.
Hence we need to prove that for any fixed $k\in\N$, 
\begin{eqnarray}
  &\biggl| \frac{d^k}{dr^k}\big|_{r=1} \det(1-r\S_n)
\biggr| \le Ce^{-c|x|^{3/2}}, \qquad & 1<\frac{2t}{n}<1+\delta, \\
  &\biggl| \frac{d^k}{dr^k}\big|_{r=1} \det(1-r\S_n)
\biggr| \le Ce^{-ct}, \qquad &1+\delta\le \frac{2t}{n}.
\end{eqnarray} 
By Cauchy's theorem, we write for $0<\epsilon<1$,
\begin{equation}
  \frac{d^k}{dr^k}\big|_{r=1} \det(1-r\S_n) 
=\frac{k!}{2\pi i} \int_{|s-1|=\epsilon}
\frac{\det(1-s\S_n)}{(s-1)^{k+1}} ds.
\end{equation} 
By Remark 3 in Section \ref{sec-iden}, for
$\varphi(z)$ as in \eqref{e5.3}, 
$\S_n$ is positive and $\|\S_n\|\le 1$. 
Hence the eigenvalues $a_j$ of (the trace class operator) 
$\S_n$ satisfies $0\le a_j\le 1$ (actually one can show 
that $0\le a_j< 1$).
For $|s-1|=\epsilon$, 
\begin{equation}
  \bigl|\det(1-s\S_n)\bigr|= \prod_j \bigl|1-sa_j\bigr|
\le \prod_j ( 1-(1-\epsilon)a_j )
= \det(1-(1-\epsilon)\S_n).
\end{equation} 
Therefore we have 
\begin{equation}
  \biggl| \frac{d^k}{dr^k}\big|_{r=1} \det(1-r\S_n) \biggr|
\le \frac{k!}{\epsilon^k} \det(1-(1-\epsilon)\S_n).
\end{equation} 
Thus it is sufficient to prove that for fixed $0<r<1$, 
\begin{eqnarray} 
  &\det(1-r\S_n)\le Ce^{-c|x|^{3/2}}, \qquad &1<\frac{2t}{n}<1+\delta,\\
  &\det(1-r\S_n)\le Ce^{-ct}, \qquad &1+\delta\le \frac{2t}{n},
\end{eqnarray} 
or by Theorem \ref{mainthm}, we need to prove that 
for any fixed $0<s<1$, (note that $\varphi$ given 
by \eqref{e5.3} has no winding)
\begin{eqnarray}
\label{eq5.9}
   &(1+s)^{-n}\det(1-s\K_n)\le Ce^{-c|x|^{3/2}}, \qquad 
&1<\frac{2t}{n}<1+\delta,\\
\label{en5.1}
   &(1+s)^{-n}\det(1-s\K_n)\le Ce^{-ct}, \qquad 
&1+\delta\le \frac{2t}{n},
\end{eqnarray} 
where $x$ is defined in \eqref{e3.4}
\begin{equation}\label{eq5.14}
  \frac{2t}{n}= 1-\frac{x}{2^{1/3}n^{2/3}}.
\end{equation}

Since $\K_n$ is an integrable operator, 
there is a naturally associated Riemann-Hilbert problem 
(see \cite{IIKS, deiftint}).
Let $m(z;k)$ be the $2\times 2$ matrix function 
which solves the following Riemann-Hilbert problem (RHP) : 
with contour $\Sigma = \{ |z|=1\}$,  
oriented counterclockwise, 
\begin{equation}\label{eq5.11}
\begin{cases}
  m(z;k) \qquad \text{is analytic in $z\in\C\setminus\Sigma$,}\\
  m_+(z;k)=m_-(z;k) \begin{pmatrix} 
1-s^2&-sz^{-k}e^{-t(z-z^{-1})}\\
sz^{k}e^{t(z-z^{-1})}&1.
\end{pmatrix} 
\qquad \text{for $z\in\Sigma$,}\\
  m(z;k) \to I \qquad \text{as $z\to\infty$.}
\end{cases}
\end{equation} 
Here and also in the following, 
the notation $f_+(z)$ (resp., $f_-(z)$) 
denotes the limiting value 
of $\lim_{z'\to z} f(z')$ from the left (resp., right) 
of the contour in the direction of the orientation.
In the above case, $m_+(z;k;s)$ (resp., $m_-(z;k;s)$) means 
$\lim_{z'\to z} m(z';k;s)$ 
with $|z'|<1$ (resp., $|z'|>1$.)
In (52) of \cite{BDJ2}, it is shown that 
\begin{equation}\label{eq5.10}
  (1+s)^{-n}\det(1-s\K_n)
= \prod_{k=n}^\infty m_{11}(0;k+1),
\end{equation}
where $m_{11}(0;k)$ is the (11)-entry of $m(z;k)$ evaluated at $z=0$.
Therefore, in order to prove \eqref{eq5.9} and \eqref{en5.1}, 
we need asymptotic results for $m_{11}(0;k)$ 
as $k,t\to\infty$.
In the special case when $s=1$, this RHP is algebraically equivalent to 
the RHP for the orthogonal polynomials on the unit circle 
with respect to the measure $e^{t(z+z^{-1})}dz/(2\pi i z)$, 
whose asymptotics as $k,t\to\infty$ was investigated in \cite{BDJ}.
The RHP \eqref{eq5.11} was introduced in \cite{BDJ2}.

There is a critical difference in the asymptotic analysis 
depending whether $s=1$ or $0<s<1$.
In the former case, the jump matrix in 
\eqref{eq5.11} has a upper/lower factorization, 
but not a lower/upper factorization, 
while in the later case, the jump matrix has both factorizations.
This difference makes the later case much easier 
to analyze asymptotically.
In the former case, we need a WKB type analysis which involves 
the construction of a parametrix in terms of the equilibrium measure 
of a certain variational problem, 
and in the case where $\frac{2t}{n}>1$, 
the main asymptotic contribution to the RHP
comes from the part of the circle near $z=-1$.
But in the later case, due to the existence of both factorization 
of the jump matrix,   
the RHP localizes in the limit just to two points on the circle.
We refer the reader to \cite{DZ1} for an example of the second type, 
and to \cite{DZ2, DVZ, DKMVZ3} for examples of the first type.

\begin{rem}
The different analysis for $s=1$ and $0<s<1$ gives us different estimates.
Indeed, when $s=1$, instead of \eqref{ee1} below, we have 
(see (6.42) of \cite{BDJ})
\begin{equation}
  \log m_{11}(0;k) \le k(-\frac{2t}{k}+\log \frac{2t}{k} +1),
\end{equation} 
which imply \eqref{q1}, \eqref{q2}.
Thus, in order to obtain the better estimates \eqref{q1}, \eqref{q2}, 
\eqref{e3.3} 
for general row $\xi_k$ (see Remark at the end of Section \ref{sec-tail}),  
we need to analysis the RHP \eqref{eq5.11} as $s\to 1$ instead of fixed $s<1$.
\end{rem}

Fix $0<s<1$.
We will prove the following estimate.
\begin{lem}\label{lem5.1}
There are positive constants $M_0$ and $c$ such that 
when $\frac{2t}{k}\ge 1+\frac{M_0}{2^{1/3}k^{2/3}}$, 
we have 
\begin{equation}\label{ee1}
  \log m_{11}(0;k) \le -c \sqrt{1-\biggl(\frac{k}{2t}\biggr)^2},
%+ C\frac{1}{k^{2/3}}.
\end{equation}
for large $t$.
\end{lem}
Assuming this result, we will prove \eqref{eq5.9} and \eqref{en5.1}, 
which completes the proof of Proposition \ref{prop3.2}.

\begin{proof}[Proof of Proposition \ref{prop3.2}]
We need to prove \eqref{eq5.9} and \eqref{en5.1}.

{\bf (a) Estimate \eqref{eq5.9} : }
It is enough to prove \eqref{eq5.9} for $-2^{1/3}n^{2/3}\delta< x\le -M$ 
where $M:=M_0(1+2\delta)$ with $M_0$ in Lemma \ref{lem5.1}. 
As noted above, since $\S_n$ is positive and $\|\S_n\|\le 1$, 
its eigenvalues $a_j$ satisfies $0\le a_j\le 1$.
Hence for $0<s<1$, by \eqref{eq5.10} and Theorem \ref{mainthm}, 
\begin{equation}
  0\le \prod_{k=n}^\infty m_{11}(0;k+1)
= \det(1-s^2\S_n) = \prod_j (1-s^2a_j) \le 1, 
\end{equation} 
for \emph{any} $n$.
Therefore in order to prove \eqref{eq5.9}, 
it is enough to show that 
\begin{equation}\label{eq5.13}
   \sum_{(*)} \log m_{11}(0;k+1) 
\le -c|x|^{3/2} + C, 
\end{equation} 
for some constants $c,C>0$ 
where the summation is over the set (note $x+M\le 0$ in \eqref{eq5.13})  
\begin{equation}\label{eq5.15}
  (*) : \quad n\le k\le n-\frac{x+M_1}{2^{1/3}}n^{1/3}, 
\end{equation} 
with $M_1:=M_0(1+\delta)$. 
We will show that for $k$ in $(*)$, 
\begin{equation}\label{eq5.16}
  1+\frac{M_0}{2^{1/3}k^{2/3}} \le \frac{2t}{k}.
\end{equation} 
Since $n\le k$, \eqref{eq5.16} follows from 
\begin{equation}\label{ee5.17}
  \frac{2t}{k} \ge 1+ \frac{M_0}{2^{1/3}n^{2/3}}.
\end{equation} 
In order to show \eqref{ee5.17}, since
\begin{equation}\label{ee5.18}
 \frac{2t}{k} = \frac{2t}{n}\cdot \frac{n}{k} 
\ge \bigl( 1-\frac{x}{2^{1/3}n^{2/3}}\bigr) 
\bigl( 1-\frac{x+M_1}{2^{1/3}n^{2/3}}\bigr)^{-1}, 
\end{equation} 
it is enough to checking that 
\begin{equation}\label{ee5.19}
 \bigl( 1-\frac{x}{2^{1/3}n^{2/3}}\bigr) 
\ge 
\bigl( 1-\frac{x+M_1}{2^{1/3}n^{2/3}}\bigr)
\biggl(  1+\frac{M_0}{2^{1/3}n^{2/3}} \biggr), 
\end{equation} 
which is equivalent to check that 
\begin{equation}\label{ee5.20}
  M_1 \ge \bigl( 1-\frac{x+M_1}{2^{1/3}n^{2/3}}\bigr) M_0.
\end{equation} 
But since $-x\le 2^{1/3}\delta n^{2/3}$, 
\begin{equation}
  \bigl( 1-\frac{x+M_1}{2^{1/3}n^{2/3}}\bigr) M_0
\le \bigl( 1+\delta -\frac{M_1}{2^{1/3}n^{2/3}}\bigr) 
\frac{M_1}{1+\delta}
\le M_1,
\end{equation} 
and hence \eqref{eq5.16} is proved.

Now using \eqref{ee1}, the sum in \eqref{eq5.13} satisfies
\begin{equation}
\begin{split}
  \sum_{(*)} \log m_{11}(0;k+1)
&\le \sum_{(*)} 
-c\sqrt{1-\biggl(\frac{k}{2t}\biggr)^2}\\
&\le \sum_{(*)} -c \sqrt{1-\frac{k}{2t}}\\
&\le -c\int_{n}^{n-\frac{x+M_1}{2^{1/3}}n^{1/3}}
\sqrt{1-\frac{s}{2t}} ds \\
&= c\biggl(\frac{4t}3\biggr) 
\biggl[ \biggr( 1-\frac{n-\frac{x+M_1}{2^{1/3}}n^{1/3}}{2t} 
\biggr)^{3/2} 
- \biggr( 1-\frac{n}{2t} \biggr)^{3/2} \biggr]\\
&= \frac23c  \biggl( \frac{2t}{n}\biggr)^{-1/2} \biggl[ 
\biggr( n^{2/3} \bigl( \frac{2t}{n}-1\bigr) + \frac{x+M_1}{2^{1/3}} 
\biggr)^{3/2} 
- \biggr( n^{2/3} \bigl( \frac{2t}{n}-1 \bigr) \biggr)^{3/2} \biggr]\\
&= \frac23c  \biggl( \frac{2t}{n}\biggr)^{-1/2} \biggl[ 
\biggr( \frac{M_1}{2^{1/3}} \biggr)^{3/2}
- \biggr( \frac{-x}{2^{1/3}} \biggr)^{3/2} \biggr]
\end{split}
\end{equation}
where the second inequality is due to the monotonicity of the function 
$f(y)=\sqrt{1-\frac{y}{2t}}$.
Since $1<\frac{2t}{n}<1+\delta$, we obtain \eqref{eq5.9}.

{\bf (b) Estimate \eqref{en5.1} : }
By a similar argument as in (a), it is enough to show that 
\begin{equation}\label{ee2}
  \sum_{(**)} \log m_{11}(0;k+1) \le -ct+C,
\end{equation}
for some constants $c,C>0$ where $(**)$ is the set 
\begin{equation}
  (**): \quad \frac{2t}{1+\delta}\le k\le \frac{2t}{1+\delta/2}.
\end{equation}
For $k$ in $(**)$, we have 
\begin{equation}
  \frac{k}{2t}
\le \frac{1}{1+\delta/2}< 1.
\end{equation}

Now using \eqref{ee1}, the sum in \eqref{ee2} satisfies
\begin{equation}
\begin{split}
  \sum_{(**)} \log m_{11}(0;k+1)
&\le \sum_{(**)} 
-c\sqrt{1-\biggl(\frac{k}{2t}\biggr)^2}\\
&\le -c \sqrt{1-\biggl(\frac{1}{1+\delta/2}\biggr)^2}
\biggl( \frac{2t}{1+\delta/2}-\frac{2t}{1+\delta}+1 \biggr) \\
& \le -ct +C.
\end{split}
\end{equation}
\end{proof}

\subsection*{RHP Asymptotics and Proof of Lemma \ref{lem5.1}}
In the rest of this section, we prove Lemma \ref{lem5.1}
by asymptotic analysis of the RHP \eqref{eq5.11}.

\begin{figure}[ht]
 \centerline{\epsfig{file=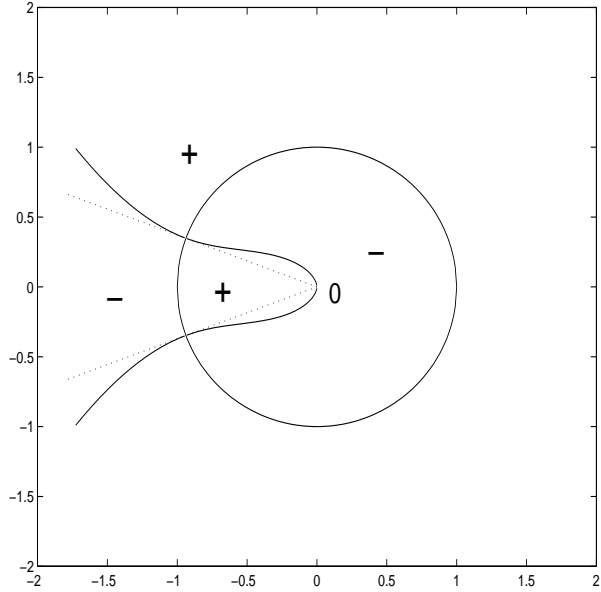, width=8cm, height=8cm}}
 \caption{Curve $Re (f(z))=0$ when $\eta=15/16$}
\label{fig-bound}
\end{figure}

Set 
\begin{equation}
  \eta := \frac{k}{2t}.
%  \gamma:= \frac{2t}{k}= 1-\frac{x}{2^{1/3}k^{2/3}}.
\end{equation} 
Under the condition of Lemma \ref{lem5.1}, 
we have $\eta<1$.
We denote by $v(z)$ the jump matrix in the second condition 
of the RHP \eqref{eq5.11}.
Note that the (21)-entry of $v$ is $se^{2tf(z;\eta)}$ where
\begin{equation}
   f(z;\eta):=\frac{1}2(z-z^{-1})+\eta\log z,
\end{equation}
where $\log z\in\R$ for $z>0$.
The critical points of this function are
$\xi:=e^{i\theta_c}$ and $\xi^{-1}=\overline{\xi}$ where 
\begin{equation}\label{er5.34}
  \xi=-\eta+i\sqrt{1-\eta^2}.
\end{equation}
Note that $-\pi/2<\theta_c<\pi$.
For $z=\rho e^{i\theta}$, consider
\begin{equation}\label{er5.35}
  F_\theta(\rho) := Re f(z) = \frac{1}2(\rho-\rho^{-1}) 
\cos\theta + \eta \log \rho.
\end{equation} 
Its derivative at $\rho=1$ is 
\begin{equation}\label{ee5.35}
  \frac{d}{d\rho}F_\theta(1)= \cos\theta+\eta, 
\end{equation} 
which is positive for $|\theta|< \theta_c$, and is negative 
for $\theta_c < |\theta|\le \pi$.
\begin{figure}[ht]
 \centerline{\epsfig{file=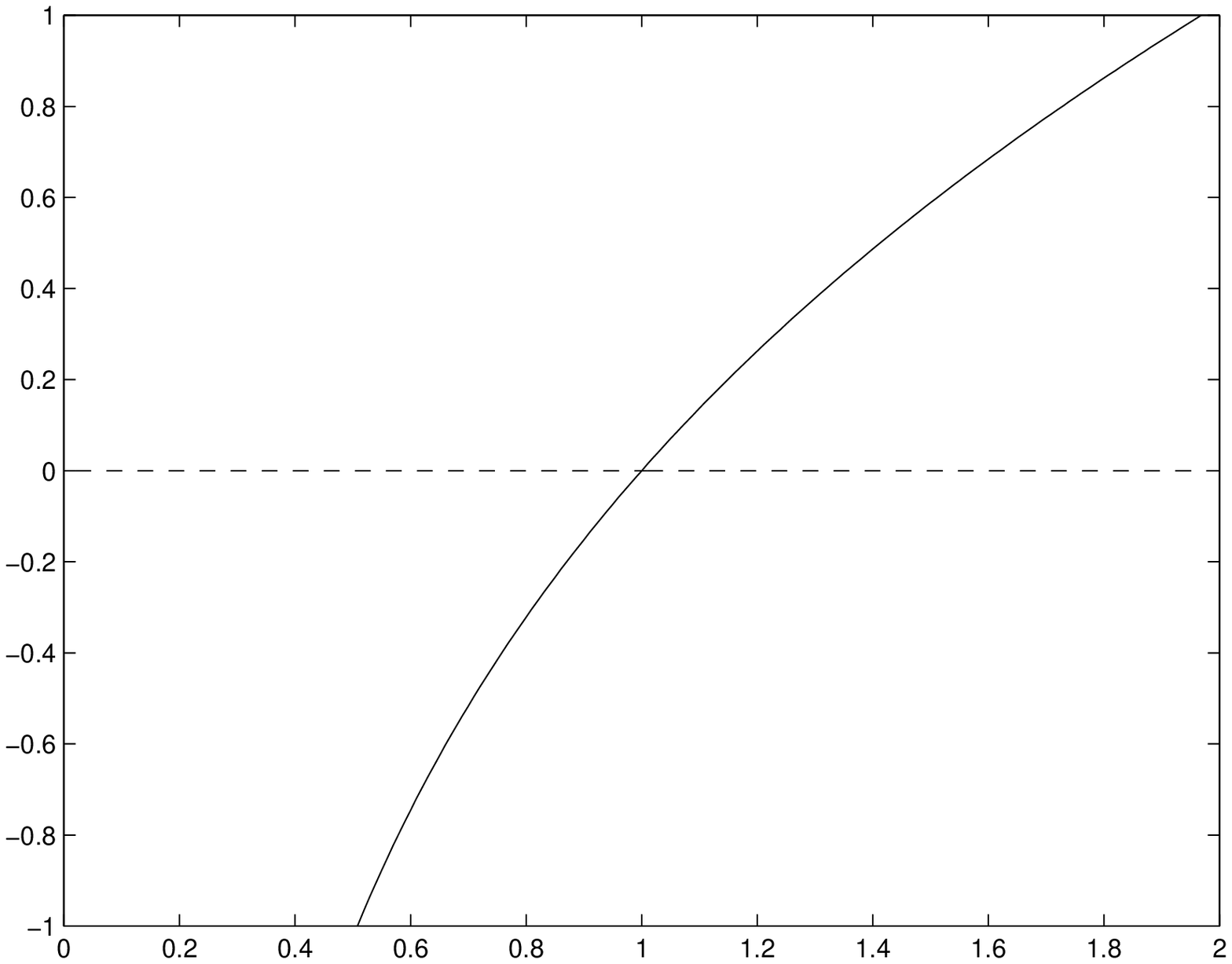, width=5cm, height=5cm}
 \epsfig{file=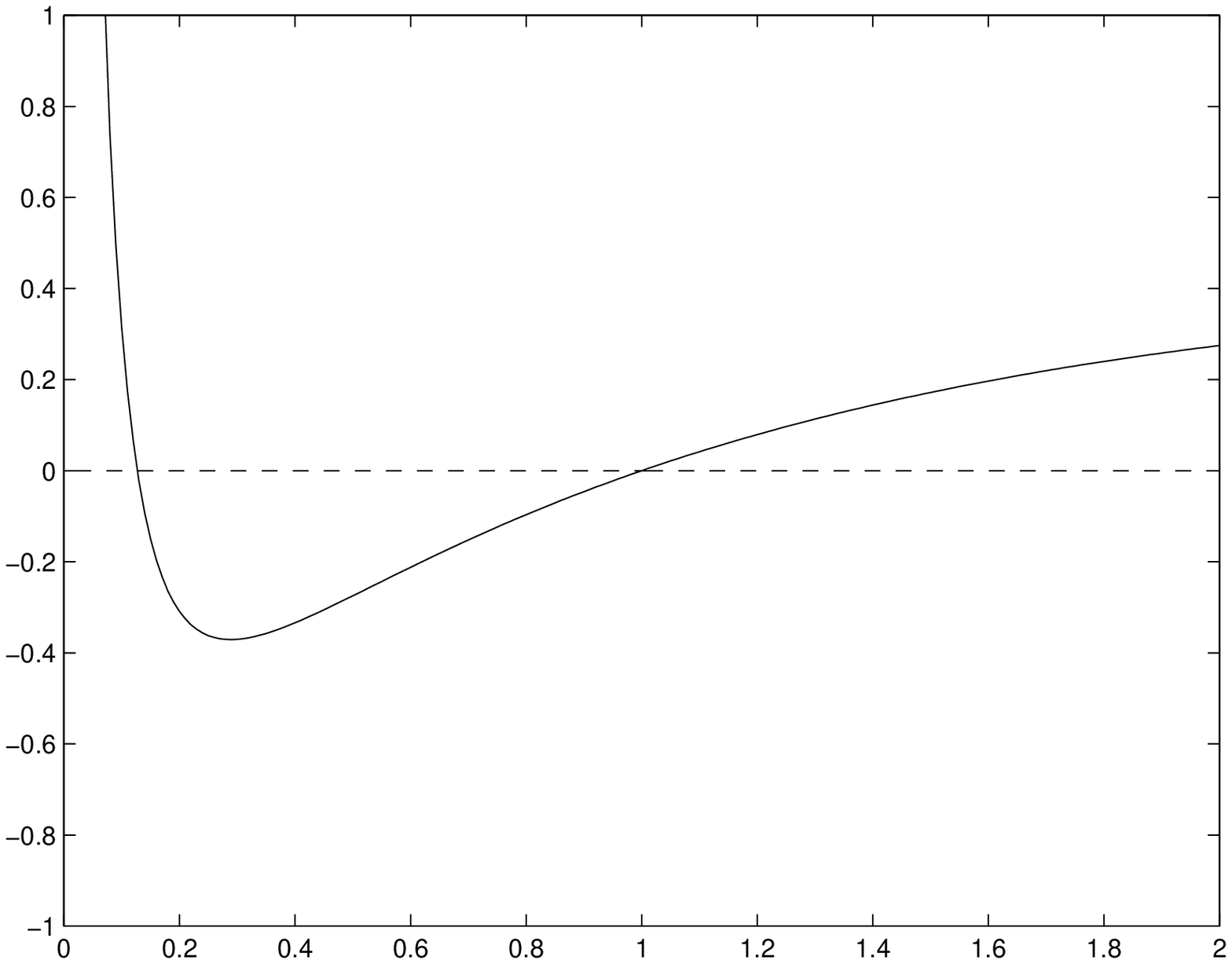, width=5cm, height=5cm}
 \epsfig{file=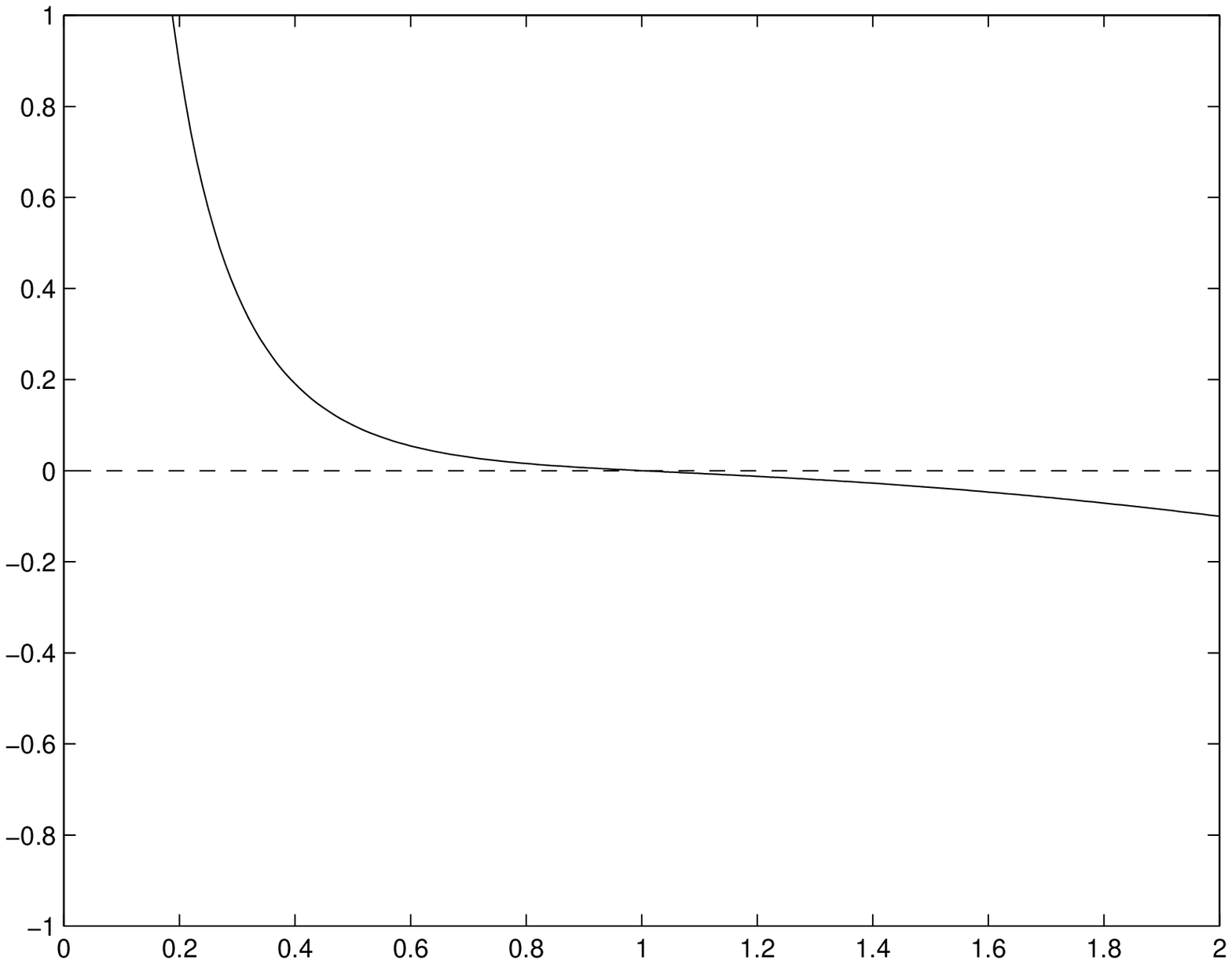, width=5cm, height=5cm}}
 \caption{Graph of $F_\theta(\rho)$ for $\theta=\frac\pi6$ (case (i)), 
$\frac{5\pi}6$ (case (ii)) and  $\pi$ (case (iii))
when $\eta=15/16$}
\label{fig-fun}
\end{figure}
Indeed it is easy to check that :
\begin{enumerate}
\item When $|\theta| \le \frac{\pi}2$, 
$F_\theta(\rho) < 0$ for $0<\rho <1$, 
and $F_\theta(\rho) > 0$ for $\rho >1$.
\item  When $\frac{\pi}2< |\theta| <\theta_c$, 
$F_\theta(\rho)>0$ for $0<\rho<\rho_0$, 
$F_\theta(\rho) <0$ for $\rho_0 <\rho <1$, 
$F_\theta(\rho) >0$ for $1<\rho <\rho_o^{-1}$, and 
$F_\theta(\rho)<0$ for $\rho>\rho_0^{-1}$. 
Here $\rho_0$ is a number satisfying 
$0<\rho_0<\rho_\theta$, where  
\begin{equation}\label{ee5.36}
  \rho_\theta:= \frac{\eta-\sqrt{\eta^2-\cos^2\theta}}{-\cos\theta} <1, 
\end{equation} 
and $\frac{d}{d\theta}F_\theta(\rho_\theta)
=\frac{d}{d\theta}F_\theta(\rho_\theta^{-1})=0$.
\item When $\theta_c < |\theta| \le \pi$, 
$F_\theta(\rho) >0$ for $0< \rho<1$,
and $F_\theta(\rho)<0$ for $\rho>1$.
\item The curve $\{ \rho e^{i\theta} : F_\theta(\rho)=0 \}$ 
crosses the circle at $90$ degree.
\end{enumerate} 
Typical graphs of $F_\theta(\rho)$ is given in Figure \ref{fig-fun} 
for the value $\eta=15/16$ for $\theta$ in the three different 
cases (i)-(iii).
Figure \ref{fig-bound} is a signature table for $Re(f(z))$ when $\eta=15/16$.
The solid curve is $Re (f(z))=0$, and the dotted rays
represent the lines $\cos\theta=-\eta$. 
The $\pm$ signs denote
the signature of $Re (f(z))$ in each of the four components.
The curve $Re (f(z))=0$ and the lines $\cos\theta=-\eta$ meet 
on the unit circle at the points $\xi$ and $\xi^{-1}$.
\begin{figure}[ht]
 \centerline{\epsfig{file=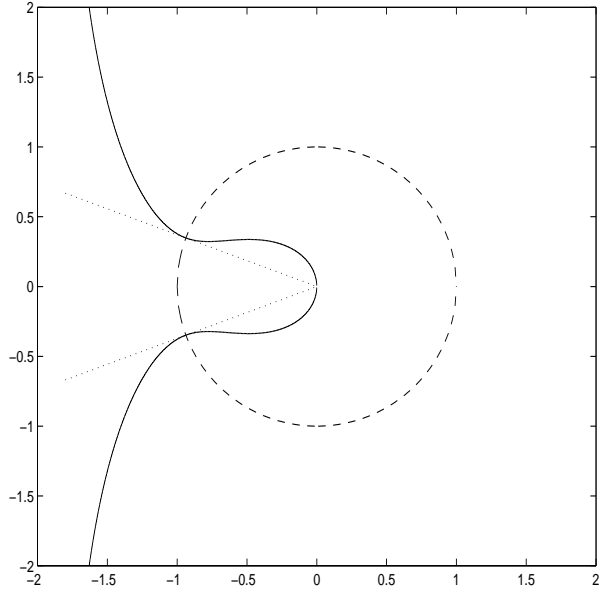, width=8cm, height=8cm}}
 \caption{Solid curve represents the set 
$\{ \rho_\theta e^{i\theta} : \frac\pi2<\theta\le \theta_c \}$ 
when $\eta=15/16$}
\label{fig-crit}
\end{figure}

Let $\Sigma=\Sigma_1\cup\Sigma_2$, 
where $\Sigma_1=\{ e^{i\theta} : |\theta| < \theta_c\}$, 
and $\Sigma_2=\Sigma\setminus\Sigma_1$.
Define the function 
\begin{equation}\label{ee5.37}
  \delta(z):=  
\biggl( \frac{z-\xi}{z-\xi^{-1}}\biggr)^{-\frac1{2\pi i} \log(1-s^2)},
\end{equation}
which is analytic in $\C\setminus\Sigma_2$ : we choose the branch so that 
$\delta(z)\to 1$ as $z\to\infty$ along the positive real axis.
Then it solves the scalar Riemann-Hilbert problem 
\begin{equation}
\begin{split} 
  &\delta_+(z)=\delta_-(z)(1-s^2), \qquad z\in\Sigma_2, \\
  & \delta(z) \to 1, \qquad \text{as $z\to\infty$,}
\end{split}
\end{equation}
where $\delta_\pm$ has the same meaning as in the RHP \eqref{eq5.11}.
Note that 
\begin{equation}
  \delta(0)= (1-s^2)^{1-\frac{\theta_c}{\pi}}. 
\end{equation}
Now set
\begin{equation}
  m^{(2)}(z):= m(z) \delta^{-\sigma_3}, 
\qquad \sigma_3= \begin{pmatrix} 1&0\\0&-1 \end{pmatrix}. 
\end{equation}
Then (i) $m^{(2)}(z)$ is analytic in $\C\setminus\Sigma$, 
(ii) $m^{(2)}(z)\to I$ as $z\to\infty$, and (iii) 
$m^{(2)}_+(z)=m^{(2)}_-v^{(2)}(z)$ for $z\in\Sigma$, where 
\begin{equation}
  v^{(2)}(z) = \begin{cases}
\begin{pmatrix} 1-s^2& -se^{-2tf(z)}\delta^2(z)\\ 
se^{2tf(z)}\delta^{-2}(z)&1 \end{pmatrix}, 
\qquad &z\in\Sigma_1, \\
\begin{pmatrix} 1& -\frac{s}{1-s^2} e^{-2tf(z)}\delta_+^2(z)\\
\frac{s}{1-s^2} e^{2tf(z)}\delta_-^{-2}&1-s^2 \end{pmatrix}, 
\qquad &z\in\Sigma_2.
\end{cases}
\end{equation}
Also we have 
\begin{equation}
  m_{11}(0)= m^{(2)}_{11}(0)(1-s^2)^{1-\frac{\theta_c}{\pi}}.
\end{equation}
Note that the jump matrix has the factorizations 
\begin{equation}
  v^{(2)}(z) = \begin{cases}
\begin{pmatrix} 1&-se^{-2tf(z)}\delta^2(z)\\0&1 \end{pmatrix}
\begin{pmatrix} 1&0\\se^{2tf(z)}\delta^{-2}(z)&1 \end{pmatrix}
, \qquad &z\in\Sigma_1, \\
\begin{pmatrix} 1&0\\
\frac{s}{1-s^2} e^{2tf(z)}\delta_-^{-2}&1-s^2 \end{pmatrix}
\begin{pmatrix} 1&-\frac{s}{1-s^2} e^{-2tf(z)}\delta_+^2(z)\\0&1 \end{pmatrix} 
, \qquad &z\in\Sigma_2.
\end{cases}
\end{equation}

\begin{figure}[ht]
 \centerline{\epsfig{file=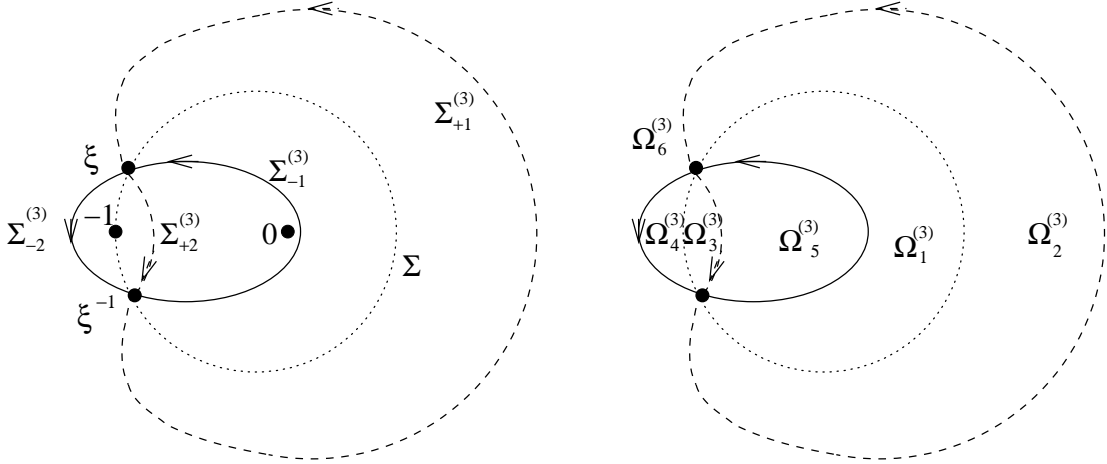, width=15cm}}
 \caption{$\Sigma^{(3)}_\pm$ and $\Omega^{(3)}_j$, $j=1,\cdots,6$}
\label{fig-contour}
\end{figure}

From (i)-(iii) in \eqref{ee5.36}, we can take an oriented closed curve
$\Sigma^{(3)}_-$
surrounding $0$ and $-1$, 
and passing through $\xi$ and $\xi^{-1}$
(the solid curve in Figure \ref{fig-contour}) on which
$Re f(z) <0$ except at $z=\xi, \xi^{-1}$.
Let $\Sigma^{(3)}_-=\overline{\Sigma^{(3)}_{-1}\cup\Sigma^{(3)}_{-2}}$
where $\Sigma^{(3)}_{-1}$ is the open subset of $\Sigma^{(3)}_-$
satisfying $|\arg(z)|<\theta_c$ and
$\Sigma^{(3)}_{-2}=\Sigma^{(3)}_-\setminus \overline{\Sigma^{(3)}_{-1}}$.
Similarly, we can take an oriented closed curve
$\Sigma^{(3)}_+$ surrounding $0$, but not $-1$
and passing through $\xi, \xi^{-1}$
(the dashed curve in Figure \ref{fig-contour}) on which
$Re f(z)>0$ except at $z=\xi, \xi^{-1}$.
Again let $\Sigma^{(3)}_+=\overline{\Sigma^{(3)}_{+1}
\cup\Sigma^{(3)}_{+2}}$
where $\Sigma^{(3)}_{+1}$ is the open subset of $\Sigma^{(3)}_+$
satisfying $|\arg(z)|<\theta_c$ and
$\Sigma^{(3)}_{+2}=\Sigma^{(3)}_+\setminus \overline{\Sigma^{(3)}_{+1}}$.
The shape of $\Sigma^{(3)}_\pm$ will be specified further below 
(see the third case for the  estimation of $|v_R(z)-I|$ 
between \eqref{en1} and \eqref{ec2}).
Let $\Sigma^{(3)}=\overline{\Sigma^{(3)}_-\cup\Sigma^{(3)}_+}$.
Let $\Omega^{(3)}_j$, $1\le j\le 6$ be open regions 
as in Figure \ref{fig-contour}.
Define 
\begin{equation}
 m^{(3)}(z) := \begin{cases}
m^{(2)}(z) \begin{pmatrix} 1&0\\ se^{2tf}\delta^{-2}&1
\end{pmatrix}^{-1}, &\qquad z\in\Omega^{(3)}_1,\\
m^{(2)}(z) \begin{pmatrix} 1&-se^{-2tf}\delta^2\\0&1
\end{pmatrix}, &\qquad z\in\Omega^{(3)}_2,\\
m^{(2)}(z) \begin{pmatrix} 1& -\frac{s}{1-s^2}e^{-2tf}\delta^2\\0&1 
\end{pmatrix}^{-1}, &\qquad z\in\Omega^{(3)}_3,\\
m^{(2)}(z) \begin{pmatrix} 1&0\\ \frac{s}{1-s^2}e^{2tf}\delta^{-2}&1
\end{pmatrix}, &\qquad z\in \Omega^{(3)}_4.\\
m^{(2)}(z), &\qquad z\in\Omega^{(3)}_5\cup\Omega^{(3)}_6.
\end{cases}
\end{equation}
Then (i) $m^{(3)}(z)$ is analytic in $\C\setminus\Sigma^{(3)}$,
(ii) $m^{(3)}(z)\to I$ as $z\to\infty$, and 
$m^{(3)}_+(z)=m^{(3)}_-(z)v^{(3)}$ for $z\in\Sigma^{(3)}$, 
where
\begin{equation}\label{e5.46}
 v^{(3)}(z)= \begin{cases}
\begin{pmatrix} 1&0\\ se^{2tf}\delta^{-2}&1
\end{pmatrix}, &\qquad z\in\Sigma^{(3)}_{-1},\\
\begin{pmatrix} 1&-se^{-2tf}\delta^2\\0&1
\end{pmatrix}, &\qquad z\in\Sigma^{(3)}_{+1},\\
\begin{pmatrix} 1& -\frac{s}{1-s^2}e^{-2tf}\delta^2\\0&1 
\end{pmatrix}, &\qquad z\in\Sigma^{(3)}_{+2},\\
\begin{pmatrix} 1&0\\ \frac{s}{1-s^2}e^{2tf}\delta^{-2}&1
\end{pmatrix}, &\qquad z\in\Sigma^{(3)}_{-2}.
\end{cases}
\end{equation}
Also we have 
\begin{equation}\label{ee5.46}
  m_{11}(0)= m^{(3)}_{11}(0)(1-s^2)^{1-\frac{\theta_c}{\pi}}.
\end{equation}
Observe that  $v^{(3)}(z)\to I$ as $t\to\infty$ 
for $z\in\Sigma^{(3)}\setminus\{\xi, \xi^{-1}\}$.
Thus we expect that $m^{(3)}(z)\to I$ as $t\to\infty$.
If this were indeed true, we would have 
\begin{equation}\label{ee5.47}
  \log m_{11}(0) \sim \frac{\sin^{-1}\sqrt{1-\eta^2}}{\pi}\log(1-s^2)
\le \frac{\log(1-s^2)}{\pi}\sqrt{1-\biggl(\frac{k}{2t}\biggr)^2}, 
\qquad t\to\infty.
\end{equation}
But the difficulty, however, is that $v^{(3)}$ does not converge to $I$ 
uniformly on $\Sigma^{(3)}$.
As in \cite{DZ1}, we overcome this difficulty by constructing a parametrix 
for the solution of the RHP $(\Sigma^{(3)},v^{(3)})$ 
around the points $\xi, \xi^{-1}$.

Let $\tau$ be a complex number satisfying $0<|\tau|<1$.
Following \cite{DZ1}, set 
\begin{equation} 
  \nu:= -\frac1{2\pi}\log(1-|\tau|^2), \qquad a:=i\nu.
\end{equation} 
Define
\begin{equation} 
  \beta_{12}:= \frac{\sqrt{2\pi}e^{\frac{\pi}{4}i}
e^{-\frac{\pi}2\nu}}{\tau\Gamma(-a)}, 
\qquad \beta_{21}:=\overline{\beta_{12}}
=\frac{\sqrt{2\pi}e^{-\frac{\pi}{4}i}
e^{-\frac{\pi}2\nu}}{\overline{\tau}\Gamma(a)}.
\end{equation} 
Note that 
\begin{equation}
  \beta_{12}\beta_{21}=\nu,
\end{equation} 
as $|\Gamma(iv)|^2=\frac{\pi}{\nu\sinh(\pi\nu)}$ for real $\nu\neq 0$.
Let $D_a$ be the parabolic-cylinder function (see, e.g. \cite{AS, WW}) 
which solves
\begin{equation}
  \frac{d^2}{d\zeta^2} D_a(\zeta)
+\biggl(\frac12-\frac{\zeta^2}{4} +a\biggr)D_a(\zeta)=0.
\end{equation} 
We note that $D_a(\zeta)$ is an entire function.

Let the matrix 
\begin{equation} 
  \Psi(w)=\begin{pmatrix} 
\Psi_{11}(w)&\Psi_{12}(w)\\ \Psi_{21}(w)&\Psi_{22}(w)
\end{pmatrix}, 
\qquad w\in\C\setminus\R,
\end{equation} 
be defined as follows (see \cite{DZ1} Section4) : for $Im(w)>0$, 
\begin{eqnarray} 
  \Psi_{11}(w) &:=& e^{-\frac34\pi\nu}D_a(e^{-\frac34\pi i}w),\\
\Psi_{12}(w) &:=&
  (\beta_{21})^{-1}e^{\frac14\pi \nu}
  \biggl(\frac{d}{dw} D_{-a}(e^{-\frac14\pi i}w)
  -\frac{iw}{2}D_{-a}(e^{-\frac14\pi i}w)\biggr), \\
\Psi_{21}(w) &:=&
  (\beta_{12})^{-1}e^{-\frac34\pi \nu}
  \biggl(\frac{d}{dw} D_{a}(e^{-\frac34\pi i}w)
  +\frac{iw}{2}D_{a}(e^{-\frac34\pi i}w)\biggr), \\
\Psi_{22}(w) &:=& e^{\frac14\pi\nu}D_{-a}(e^{-\frac14\pi i}w),
\end{eqnarray} 
and for $Im(w)<0$, 
\begin{eqnarray} 
  \Psi_{11}(w) &:=& e^{\frac14\pi\nu}D_a(e^{\frac14\pi i}w),\\ 
\Psi_{12}(w) &:=&
  (\beta_{21})^{-1}e^{-\frac34\pi \nu}
  \biggl(\frac{d}{dw} D_{-a}(e^{\frac34\pi i}w)
  -\frac{iw}{2}D_{-a}(e^{\frac34\pi i}w)\biggr), \\
\Psi_{21}(w) &:=&
  (\beta_{12})^{-1}e^{\frac14\pi \nu}
  \biggl(\frac{d}{dw} D_{a}(e^{\frac14\pi i}w)
  +\frac{iw}{2}D_{a}(e^{\frac14\pi i}w)\biggr), \\
\Psi_{22}(w) &:=& e^{-\frac34\pi\nu}D_{-a}(e^{\frac34\pi i}w),
\end{eqnarray} 
The function $\Psi$ satisfies
\begin{itemize}
\item $\Psi(w)$ is analytic in $w\in\C\setminus\R$.
\item For $w\in\R$, 
\begin{equation}
  \Psi_+(w)= \Psi_-(w) \begin{pmatrix} 
1-|\tau|^2&-\overline{\tau}\\ \tau&1
\end{pmatrix}, 
\end{equation} 
where $\Psi_+(w)$ (resp., $\Psi_-(w)$) 
is the limit of $\Psi(s)$ as $s\to w$ with $Im(s)>0$ (resp., $Im(s)<0$).
\item As $w\to\infty$, 
\begin{equation}\label{ea5.63}
  \Psi(w)e^{\frac14iw^2\sigma_3}w^{-i\nu\sigma_3} = I + O(w^{-1})
\end{equation} 
where $w^{-i\nu}$ denotes the branch which 
is analytic in $\C\setminus(-\infty,0]$ 
and has modulus $1$ for $w\in (0,\infty)$.
\end{itemize}
These properties can be found in \cite{DZ1} Section 4.

\begin{figure}[ht]
 \centerline{\epsfig{file=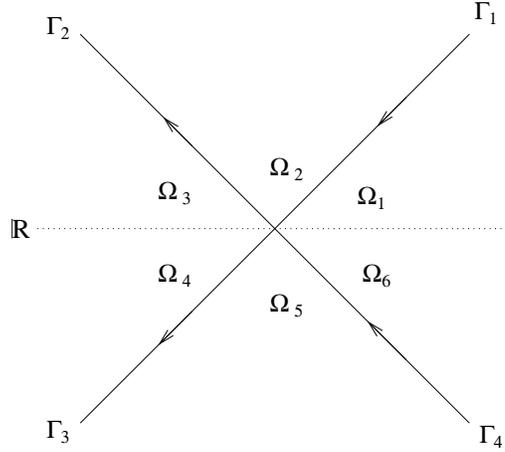, width=7cm}}
 \caption{$\Gamma_j$ and $\Omega_j$}
\label{fig-Gamma}
\end{figure}

Let $\Gamma$ be the union of four rays, labeled by 
$\Gamma_j$, $j=1,\cdots,4$, 
with the orientation as indicated in 
Figure \ref{fig-Gamma}.
All the rays and $\R$ meet at the angle $\pi/3$.
Denote the components of $\C\setminus(\Gamma\cup\R)$ by $\Omega_j$, 
$j=1,\cdots,6$ as in Figure \ref{fig-Gamma}.
Define $H(w)$, analytic in $\C\setminus\Gamma$, by 
\begin{equation}
  H(w) := \bigl(\begin{smallmatrix} a(\xi)^{-1}&0\\0&a(\xi) 
\end{smallmatrix}\bigr)
\bigl(\begin{smallmatrix} 0&1\\1&0 \end{smallmatrix}\bigr)
 \Psi(w) e^{\frac14 iw^2\sigma_3} w^{-i\nu\sigma_3}
\bigl(\begin{smallmatrix} 0&1\\1&0 \end{smallmatrix}\bigr)
\phi(w)
\bigl(\begin{smallmatrix} a(\xi)&0\\0&a(\xi)^{-1} 
\end{smallmatrix}\bigr),
\end{equation} 
where $a(\xi)$ is 
\begin{equation}
  a(\xi)=e^{tf(\xi)} \biggl( 
\frac{-i\xi}{(\xi-\xi^{-1})\sqrt{2t}(1-\eta^2)^{1/4}} \biggr)^{i\nu},
\end{equation} 
and $\phi(w)$ is defined by 
\begin{equation}
\phi(w):= 
\begin{cases} 
   \begin{pmatrix} 1&-\tau e^{\frac12 iw^2} w^{-2i\nu}\\ 0&1 
\end{pmatrix} , \qquad
&w\in \Omega_1, \\
   \begin{pmatrix} 1&0\\ 
-\overline{\tau}e^{-\frac12 iw^2} w^{2i\nu}&1 
\end{pmatrix} , \qquad 
&w\in \Omega_6, \\
   \begin{pmatrix} 1&0\\ \frac{\overline{\tau}}{1-|\tau|^2}
e^{-\frac12 iw^2} w^{2i\nu}&1 \end{pmatrix} ,
\qquad &w\in \Omega_3, \\
   \begin{pmatrix} 1&\frac{\tau}{1-|\tau|^2}
e^{\frac12 iw^2} w^{-2i\nu} \\ 0&1 
\end{pmatrix} , \qquad 
&w\in \Omega_4, \\
   I
&w\in \Omega_2, \Omega_5.
\end{cases} 
\end{equation} 
Then by recalling that $w^{-i\nu\sigma_3}$ is analytic in 
$\C\setminus(-\infty,0]$, one can directly check that 
$H_+(w)=H_-(w)v_H(w)$ for $w\in\Gamma$, where $v_H(w)$ is given by 
\begin{equation}\label{e5.67}
v_H(w):=
\begin{cases}
   \begin{pmatrix} 1&-\tau a(\xi)^{-2} w^{-2i\nu}e^{\frac12 iw^2}
\\ 0&1 \end{pmatrix} , \qquad
&w\in \Gamma_1, \\
   \begin{pmatrix} 1&0\\ \overline{\tau} a(\xi)^{2} w^{2i\nu}e^{-\frac12 iw^2}
&1 \end{pmatrix} , \qquad
&w\in \Gamma_4, \\
   \begin{pmatrix} 1&0\\ 
\frac{\overline{\tau}}{1-|\tau|^2}a(\xi)^{2} w^{2i\nu}e^{-\frac12 iw^2}
&1 \end{pmatrix} ,
\qquad &w\in \Gamma_2, \\
   \begin{pmatrix} 1&\frac{-\tau}{1-|\tau|^2} a(\xi)^{-2} 
w^{-2i\nu}e^{\frac12 iw^2} \\ 0&1
\end{pmatrix} , \qquad
&w\in \Gamma_3.
\end{cases}
\end{equation}
Also, from \eqref{ea5.63}, we have 
\begin{equation}\label{ee5.67}
  H(w) = I + O(w^{-1}), \qquad \text{as $w\to\infty$.}
\end{equation} 
As $|a(\xi)|=e^{-\nu\theta_c}$, $-\pi/2<\theta_c<\pi$, we see that the error 
term $O(w^{-1})$ in \eqref{ee5.67} is uniform for $\frac{2t}{k}>1$.
Similarly, $|H(w)|$ is uniformly bounded in the $w$ plane for 
$\frac{2t}{k}>1$.

\begin{figure}[ht]
 \centerline{\epsfig{file=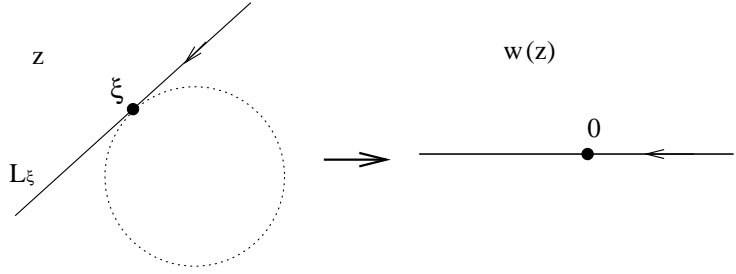, width=10cm}}
 \caption{map $z\mapsto w(z)$}
\label{fig-map}
\end{figure}
Define the map 
\begin{equation}
   z \mapsto w(z):=\sqrt{2t}(1-\eta^2)^{1/4}i\xi^{-1}(z-\xi).
\end{equation} 
It maps $\xi$ to $0$, and the tangent line $L_\xi$ to the 
unit circle $\Sigma$ at $\xi$, to the real line as in Figure \ref{fig-map}.
Let $\calO_\xi$, $\calO_{\overline{\xi}}$  
be the disjoint sets $\{z : |z-\xi|<\varrho \}$, 
$\{z : |z-\overline{\xi} |<\varrho \}$, respectively,  
where $\varrho$ 
is defined by 
\begin{equation}\label{ee5.69}
  \varrho= \begin{cases} 
\frac{\epsilon}{2}|\xi-\overline{\xi}| = \epsilon\sqrt{1-\eta^2}, \qquad 
& 1+\frac{M_0}{2^{1/3}k^{2/3}}\le \frac{2t}{k}< 1+\delta,\\
\epsilon, \qquad & 1+\delta \le \frac{2t}{k}.
\end{cases} 
\end{equation} 
The (small) parameter $0<\epsilon<1$ will be specified below (see 
\eqref{ee5.97} below).
We note that one may choose the curves in $\Sigma^{(3)}$ above so that 
in $\calO_\xi$, $\calO_{\overline{\xi}}$, they 
are straight lines 
which map under $z\mapsto w(z)$ to (finite subsets of ) 
the rays $\Gamma_j$, $j=1,\cdots, 4$,
$\Sigma^{(3)}_{-1}\cap\calO_\xi \to \Gamma_4$, 
$\Sigma^{(3)}_{-2}\cap\calO_\xi \to \Gamma_2$, 
$\Sigma^{(3)}_{+1}\cap\calO_\xi \to \Gamma_1$, 
$\Sigma^{(3)}_{+2}\cap\calO_\xi \to \Gamma_3$, 
and similarly for the neighborhood of $\calO_{\overline{\xi}}$.
For $\tau=s$, we define 
\begin{equation}
  m_p(z) := \begin{cases} 
H(w(z)), \qquad &z\in \calO_\xi\setminus \Sigma^{(3)}, \\
\overline{H(w(\overline{z}))}, 
&z\in \calO_{\overline{\xi}}\setminus \Sigma^{(3)}, \\
I, &z\in \C\setminus 
\overline{(\calO_\xi\cup \calO_{\overline{\xi}})}. \\
\end{cases} 
\end{equation} 
\begin{figure}[ht]
 \centerline{\epsfig{file=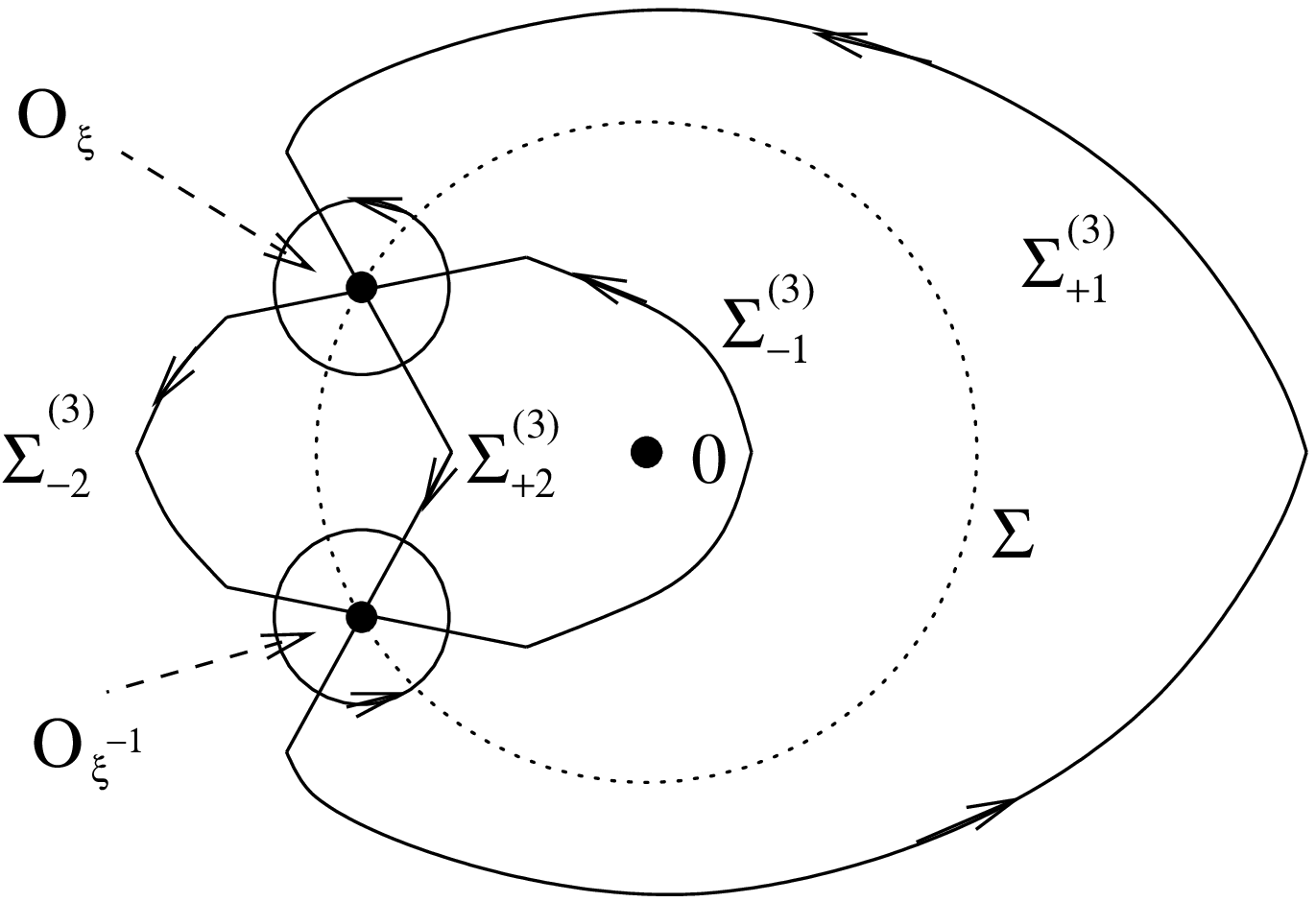, width=8cm}}
 \caption{$\Sigma_R:= \Sigma^{(3)}\cup\partial \calO_\xi\cup\partial
\calO_{\overline{\xi}}$}
\label{fig-sigmar}
\end{figure}
Let $\Sigma_R:= \Sigma^{(3)}\cup\partial \calO_\xi\cup\partial 
\calO_{\overline{\xi}}$ as in Figure \ref{fig-sigmar}
where $\partial \calO_\xi$ and $\calO_{\overline{\xi}}$ are
oriented counterclockwise.
Clearly $m_p$ solves a RHP on $\Sigma_R$ : $m_p(z)$ is 
analytic in $\C\setminus \Sigma^{(3)}$, $m_p(z)\to I$ as $z\to\infty$, 
and $m_{p+}(z)=m_{p-}(z)v_p(z)$ for $z\in\Sigma_R$ 
for a suitable jump matrix $v_p$.
Set $R(z):= m^{(3)}(z)m_p(z)^{-1}$.
Then $R_+(z)=R_-(z)v_R(z)$ for $z\in\Sigma_R$ where 
$v_R=m_{p-}v^{(3)}v_p^{-1}m_{p-}^{-1}$.
Now we estimate $|v_R(z)-I|$.

\begin{itemize}                  %%  main estimates %%
\item 
For $z\in\partial \calO_\xi$, we have 
from \eqref{ee5.69}, 
when $1+\frac{M_0}{2^{1/3}k^{2/3}}\le \frac{2t}{k}\le 1+\delta$,
\begin{equation}
\begin{split} 
  |w(z)| &= \epsilon(1-\eta^2)^{3/4}\sqrt{2t} \\
&\ge \epsilon  \biggl(\frac{2t}{k}\biggr)^{-1/4} 
\biggl( k^{2/3}\bigl(\frac{2t}{k}-1\bigr)\biggr)^{3/4} \\
&\ge \frac{\epsilon M_0^{3/4}}{2^{1/4}(1+\delta)^{1/4}}.
\end{split} 
\end{equation} 
When $1+\delta\le \frac{2t}{k}$, 
\begin{equation}
  |w(z)| = \sqrt{2t}(1-\eta^2)\epsilon
\ge \biggl( \frac{\delta}{1+\delta}\biggr)^{1/4}\epsilon \sqrt{2t}.
\end{equation} 
Thus if we have taken $M_0$ large, and $t$ is large, we have 
for $z\in \partial \calO_\xi$, 
from \eqref{ee5.67}, 
\begin{equation} 
  m_{p+}(z) = I+ O\biggl(\frac{1}{\min(M_0^{3/4},\sqrt{t})}\biggr).
\end{equation} 
But as $v^{(3)}=I$ on $\partial\calO_\xi$, 
$v_R(z)=v_p(z)^{-1}= m_{p+}(z)^{-1}$, and hence 
\begin{equation}\label{ee5.72}
  \|v_R(z)-I\|_{L^\infty(\partial\calO_\xi)}
\le \frac{C}{\min(M_0^{3/4},\sqrt{t})}
\qquad z\in \partial \calO_\xi.
\end{equation} 
We are using here the standard fact that if $\det v_p=1$, then $\det m_p=1$.
Similarly, we have the same estimate \eqref{ee5.72} on $v_R(z)$ for 
$z\in \partial \calO_{\overline{\xi}}$. 
\item For $z\in \Sigma^{(3)}\cap \calO_\xi$, 
since $m_p$ and $m_p^{-1}$ are uniformly bounded, 
$|v_R(z)-I| \le C|v^{(3)}(z)v_p(z)^{-1} -I|$.
For $z\in \Sigma^{(3)}_{-1} \cap \calO_\xi$, 
by \eqref{e5.46}, \eqref{e5.67}, 
\begin{equation}\label{ee5.74}
\begin{split} 
  |v^{(3)}(z)v_p(z)^{-1} -I|
&\le | v^{(3)}_{21}(z)- (v_H)_{21}(w(z))|\\
&= s| e^{2tf(z)}\delta^{-2}(z)-a(\xi)^2w(z)^{2i\nu}e^{-\frac12iw(z)^2}|
=:s |\Delta|.
\end{split} 
\end{equation} 
Setting $u:=i\xi^{-1}(z-\xi)$, we have 
\begin{equation}
\begin{split} 
  \Delta &= e^{2t(\frac12(\xi(1-iu)-\frac1{\xi(1-iu)})+\eta\log \xi(1-iu))}
\biggl( \frac{-i\xi u}{\xi(1-iu)-\xi^{-1}} \biggr)^{2i\nu} \\
&\qquad - e^{2t(\frac12(\xi-\xi^{-1})+\eta\log \xi)-ti(1-\eta^2)^{1/2}u^2}
\biggl( \frac{-i\xi u}{\xi-\xi^{-1}} \biggr)^{2i\nu} \\
&=\bigl(  e^{2th(u)}j(u)-1 \bigr) 
e^{2t(\frac12(\xi-\xi^{-1})+\eta\log \xi)-ti(1-\eta^2)^{1/2}u^2}
\biggl( \frac{-i\xi u}{\xi-\xi^{-1}} \biggr)^{2i\nu},
\end{split} 
\end{equation} 
where 
\begin{equation}
\begin{split} 
  h(u) &= \frac12 ( -i\xi u -\frac{i u}{\xi(1-iu)} ) 
+ \eta\log(1-iu) + \frac12 i(1-\eta^2)^{1/2}u^2  \\
&= (-\frac16 i\eta+\frac12\sqrt{1-\eta^2}) u^3 +O(u^4),
\end{split} 
\end{equation} 
and 
\begin{equation}
  j(u)= \biggl( \frac{\xi-\xi^{-1}}{\xi-\xi^{-1}-i\xi u} \biggr)^{2i\nu}
= 1+ O(\frac{u}{\xi-\xi^{-1}}) = 1+O(\frac{u}{\sqrt{1-\eta^2}}).
\end{equation} 
Also, as $|\frac{u}{\sqrt{1-\eta^2}}| \le c\epsilon$ for $z\in \calO_\xi$, 
we have 
\begin{eqnarray} 
  | h(u)| &\le & c |u|^3, \\
%\frac12\sqrt{1-\eta^2} |u|^3+ O(|u|^4)
%= c\sqrt{1-\eta^2}|u|^3, \\
  j(u) &=& 1+O(\frac{u}{\sqrt{1-\eta^2}}) = 1+O(\epsilon).
\end{eqnarray} 
On the other hand, for $z\in\Sigma^{(3)}_{-1}\cap\calO_\xi$, 
\begin{equation}
   Re (-it\sqrt{1-\eta^2} u^2) \le -ct\sqrt{1-\eta^2} |u|^2,
\qquad c=\cos\frac\pi6>0.
\end{equation} 
Therefore, we obtain 
\begin{equation}
\begin{split} 
  |\Delta| 
&\le C|(e^{2th(u)}-1)+e^{2th(u)}(j(u)-1)| e^{-ct\sqrt{1-\eta^2} |u|^2}\\
&\le C\bigl(|2th(u)| + |j(u)-1|\bigr)
e^{-ct\sqrt{1-\eta^2} |u|^2 + |Re (2th(u))|} \\
&\le C(t|u|^3 + \epsilon) e^{-ct\sqrt{1-\eta^2} |u|^2}\\
&\le \frac{C}{\sqrt{2t}(1-\eta^2)^{3/4}} + C\epsilon,
\end{split} 
\end{equation} 
where for the last inequality, we have used the fact that 
$|x^3e^{-x^2}|$ is uniformly bounded for $x\in\R$.
Now for $1+\frac{M_0}{2^{1/3}k^{2/3}} \le \frac{2t}{k} <1+\delta$, 
we have 
\begin{equation}
  \sqrt{2t}(1-\eta^2)^{3/4} 
> \sqrt{2t}(1-\eta)^{3/4}
= \biggl( \frac{k}{2t} \biggr)^{1/4} 
\biggr( k^{2/3}(\frac{2t}{k}-1)\biggr)^{3/4} \ge 
  \frac{1}{(1+\delta)^{1/4}}\biggl( \frac{M_0}{2^{1/3}} \biggr)^{3/4}, 
\end{equation} 
and for $1+\delta \le \frac{2t}{k}$, we have 
\begin{equation}
  \sqrt{2t}(1-\eta^2)^{3/4} \ge \sqrt{2t} \biggl( \frac{\delta}{1+\delta}
\biggr)^{3/4}.
\end{equation} 
Thus we obtain (recall \eqref{ee5.74})
\begin{equation}\label{er5.86}
  \|v_R-I\|_{L^{\infty}(\Sigma^{(3)}_{-1}\cap\calO_\xi)}
\le Cs\|\Delta\|_{L^{\infty}(\Sigma^{(3)}_{-1}\cap\calO_\xi)} 
\le \frac{C}{\min (M_0^{3/4}, \sqrt{t})} + C\epsilon,
\end{equation} 
which is small if we take $M_0, t$ large and $\epsilon$ small.
For other parts of $\Sigma^{(3)}\cap\calO_\xi$, by a similar argument, 
we obtain the same estimate.
By the symmetry $m^{(3)}(z)=\overline{m^{(3)}(\overline{z})}$ 
and $m_p(z)=\overline{m_p(\overline{z})}$, we obtain the same estimate 
for $\Sigma^{(3)}\cap\calO_{\overline{\xi}}$.
\item 
Let $\calO:=\calO_\xi\cup\calO_{\overline{\xi}}$.
For $z\in \Sigma^{(3)}_{-1}\cap(\C\setminus\overline\calO)$,
$v_R(z)=v^{(3)}(z)$. 
Thus we need an estimate for $v^{(3)}_{21}(z)=se^{2tf(z)}\delta^{-2}(z)$.
Since $|\delta(z)|= e^{-\nu \theta}$ 
where $\theta=arg (\frac{z-\xi}{z-\overline{\xi}})$, 
$|\delta(z)|$ and $|\delta^{-1}(z)|$ are uniformly bounded.
When $1+\delta\le \frac{2t}{k}$, 
$dist(\Sigma^{(3)}_{-1}\cap\overline{\calO}^{c}, \{ \xi, \xi^{-1}\}$ 
is uniformly bounded below.
From this fact, one can check that we can take $\Sigma^{(3)}_{-1}$ 
so that $Re (f(z))\le -c_0(\epsilon)$ 
for $z\in\Sigma^{(3)}_{-1}\cap\overline{\calO}^{c}$
for some constant $c_0(\epsilon)>0$ depending on $\epsilon$.
Hence we have 
\begin{equation}\label{en1}
   \| v_R -I \|_{L^\infty(\Sigma^{(3)}_{-1}\cap\overline{\calO}^{c})}
\le Ce^{-c_0(\epsilon) t}, 
\qquad 1+\delta\le \frac{2t}{k}.
\end{equation} 
On the other hand, when $1+\frac{M_0}{2^{1/3}k^{2/3}}\le 
\frac{2t}{k}<1+\delta$, we take 
$\Sigma^{(3)}_{-1}=\{ \rho(\theta) e^{i\theta} : |\theta|< \theta_c\}$ 
such that 
\begin{enumerate} 
\item[(i)] For $z\in\Sigma^{(3)}_{-1}$ with 
$\frac{2\pi}3< | arg(z)| < \theta_c$, $\Sigma^{(3)}_{-1}$ is 
a pair of straight 
lines which meet the unit circle at $\xi$ and $\overline{\xi}$, respectively, 
with angle $\pi/3$.
\item[(ii)] For $\rho e^{i\theta}\in \Sigma^{(3)}_{-1}$ 
with $| arg(z)|\le \frac{2\pi}3$, 
$Re (f(\rho e^{i\theta})) \le Re (f(\rho' e^{i\theta}))$ 
for $\rho\le \rho'\le 1$.
Also $\rho(\theta)$ is an increasing function for $0<\theta<\frac{2\pi}3$ 
and is a decreasing function for $-\frac{2\pi}3<\theta<0$.
\end{enumerate} 
(Here the precise value $\frac{2\pi}3$ is of no importance : any angle 
between $\pi/2$ and $\pi$ will do.) 
Condition (ii) can be achieved by choosing $\Sigma^{(3)}_{-1}$ 
always to be above the curve $\{ \rho_\theta e^{i\theta} : \frac{\pi}2 < 
\theta \le \theta_c \}$ (recall \eqref{ee5.36} and Figures 
\ref{fig-fun}, \ref{fig-crit}).
Condition (i) can be achieved as the curve 
$\{ \rho_\theta e^{i\theta} : \frac{\pi}2 < \theta \le \theta_c \}$
crosses the unit circle at $90$ degree (see Figure \ref{fig-crit})
For $z$ in (i) satisfying $arg(z)>0$, we have 
$z=\xi(1-i r e^{-\frac{\pi}3 i})$ for some real $r>0$.
We note that $r\le \frac{\sqrt{3}}2<\frac{2}{\sqrt{3}}$.
For such $z$, we have (recall \eqref{er5.34})
\begin{equation}\label{ec2}
  Re (f(z))= A(r)\sqrt{1-\eta^2} + B(r) \eta, 
\end{equation} 
where
\begin{equation} 
A(r)= \frac{r^2(r-\sqrt3)}{4(1-\sqrt3 r +r^2)}, \qquad 
B(r)= -\frac{r(r-\sqrt3)(2-\sqrt3 r)}{4(1-\sqrt3 r +r^2)}
+ \frac12 \log(1-\sqrt3 r +r^2).
\end{equation} 
One can easily check that $A(r)<0$ for $0<r<\sqrt{3}$ 
and $B(r)<0$ for $0<r<\frac{2}{\sqrt{3}}$.
Thus for $z$ in (i) satisfying $arg(z)>0$, we have 
for some $c>0$, 
\begin{equation}\label{ez2}
  Re (f(z)) \le A(r)\sqrt{1-\eta^2}\le -cr^2\sqrt{1-\eta^2}.
\end{equation} 
For $z$ in (ii), note first that  
for fixed $0<\rho<1$, $Re (f(\rho e^{i\theta})$ is an increasing function
in $0\le \theta<\pi$.
Let $z_b$ be the point on $\Sigma^{(3)}_{-1}$ satisfying $arg(z)=\frac{2\pi}3$. 
Thus together with the condition (ii), 
we obtain for $z$ in (ii) satisfying $arg(z)>0$, 
\begin{equation}
  Re (f(z)) \le Re (f(z_b)) \le -c|z_b-\xi|^2\sqrt{1-\eta^2} 
\le -c|z-\xi|^2\sqrt{1-\eta^2}.
\end{equation} 
Here the second inequality follows from \eqref{ez2}.
Thus we have 
for $z \in \Sigma^{(3)}_{-1}\cap\overline{\calO_\xi}^c$
with $arg(z)>0$, 
\begin{equation}\label{er92}
  |v_R(z)-I| \le 
\begin{cases} 
Ce^{-c_0(\epsilon) t}, \qquad &1+\delta\le \frac{2t}{k},\\
Ce^{-ct|z-\xi|^2\sqrt{1-\eta^2}}, \qquad &1+\frac{M_0}{2^{1/3}k^{2/3}}
\le \frac{2t}{k}< 1+\delta. 
\end{cases} 
\end{equation} 
By symmetry, we have similar estimates for 
$z\in \Sigma^{(3)}_{-1}$ with $arg(z)<0$.
Since $|z-\xi|> \varrho$ for $z\in\overline{\calO}^c$, 
the above estimates imply in particular that 
\begin{equation}\label{er5.89}
  \| v_R -I\|_{L^{\infty}(\Sigma^{(3)}_{-1}\cap(\C\setminus\overline\calO))}
\le 
\begin{cases}
Ce^{-c\epsilon^2 M_0^{3/2}}, \qquad
& 1+\frac{M_0}{2^{1/3}k^{2/3}}\le \frac{2t}{k}<1+\delta, \\
Ce^{-c_0(\epsilon) t}, \qquad & 1+\delta\le \frac{2t}{k},
\end{cases}
\end{equation}
For $\Sigma^{(3)}_{+1}\cap(\C\setminus\overline\calO)$,
by the symmetry $Re (f(\rho e^{i\theta})) = Re (f(\rho^{-1} e^{i\theta})$,
we have the same estimate.
Also by a similar argument, we obtain a similar estimate
for $(\Sigma^{(3)}_{-2}\cup\Sigma^{(3)}_{+2})\cap(\C\setminus\overline\calO)$.
\end{itemize}

As usual, define an operator on $L^2(\Sigma_R)$,
\begin{equation}
  C_{v_R}(f) = C_-(f(v_R-I))
\end{equation}
where $C_-$ is the Cauchy operator
\begin{equation}
  (C_-f)(z)= \lim_{z'\to z} \frac1{2\pi i}\int_{\Sigma_R} \frac{f(s)}{s-z'} ds,
\qquad \text{$z\in\Sigma_R$,\ \ $z'$ on the $-$ side of $\Sigma_R$.}
\end{equation}
As the Cauchy operator is scale invariant, 
$C_-$ is bounded from $L^2(\Sigma_R)\to L^2(\Sigma_R)$ 
uniformly for $\frac{2t}{k}\ge 1+\frac{M_0}{2^{1/3}k^{2/3}}$,
%and $\|v_R-I\|_{L^\infty(\Sigma_R)}\le Ct^{-1/2}+Ce^{-ct^{1/3}}<1$
%when $t$ is large, 
and we have $\|C_{v_R}\|<\frac12$ 
for $t$, $M_0$ sufficiently large by \eqref{ee5.72}, \eqref{er5.86} 
and \eqref{er5.89}.
Hence $1-C_{v_R}$ is invertible.
By standard facts in Riemann-Hilbert theory (see \cite{CG, BC}),
the solution $R(z)$ to the RHP $(\Sigma_R,v_R)$ is given by
\begin{equation}
  R(z)= I+ \frac1{2\pi i}
\int_{\Sigma_R} \frac{(I+(1-C_{v_R})^{-1}C_{v_R}I)(v_R-I)(s)}{s-z} ds.
\end{equation}
As $m_p(0)=I$, we have $m^{(3)}_{11}(0)=R_{11}(0)$.
By using $dist (0, \Sigma_R)>0$,
$\|(1-C_{v_R})^{-1}\|\le c$,
and $\|C_-\|\le c$, we have
\begin{equation}\label{ee5.92}
\begin{split}
 |m^{(3)}_{11}(0)-1|
&\le c\|v_R-I\|_{L^1(\Sigma_R)}
+ c\|(1-C_{v_R})^{-1}C_{v_R}I\|_{L^2(\Sigma_R)}\|v_R-I\|_{L^2(\Sigma_R)}\\
&\le c\|v_R-I\|_{L^1}
+ c\|(1-C_{v_R})^{-1}\|_{L^2\to L^2}\|C_-(v_R-I)\|_{L^2}
\|v_R-I\|_{L^2}\\
&\le c\|v_R-I\|_{L^1}
+ c \|v_R-I\|^2_{L^2}\\
&\le c\|v_R-I\|_{L^1}
+ c \|v_R-I\|_{L^\infty} \|v_R-I\|_{L^1}\\
&\le c\|v_R-I\|_{L^1(\Sigma_R)}
\end{split}
\end{equation}
as $\|v_R-I\|_{L^\infty}$ is bounded.
We estimate $\|v_R-I\|_{L^1}$ in each part of $\Sigma_R$.
First, for $\partial\calO$ and $\Sigma^{(3)}\cap\calO$,
since the length of the contour is of order $\varrho$, we obtain
by \eqref{ee5.72}, \eqref{er5.86}
%(we take $\sqrt{t}>M_0^{3/4}$)
\begin{equation}\label{ee5.93}
  \|v_R-I\|_{L^1(\Sigma_R\cap\overline{\calO})}
\le C\varrho\biggl(\frac1{\min (M_0^{3/4}, \sqrt{t})}+\epsilon\biggr). 
\end{equation}
When $1+\frac{M_0}{2^{1/3}k^{2/3}}\le \frac{2t}{k}<1+\delta$, 
by \eqref{ee5.69}, $\varrho= \epsilon \sqrt{1-\eta^2}$.
When $1+\delta \le \frac{2t}{k}$, 
$\sqrt{1-\eta^2}=\sqrt{1-(k/(2t))^2} \ge C$, 
and hence we have $\varrho=\epsilon\le c\sqrt{1-\eta^2}$.
Thus in both cases, we obtain 
\begin{equation}\label{ee5.94}
  \|v_R-I\|_{L^1(\Sigma_R\cap\overline{\calO})}
\le C\sqrt{1-\eta^2}\biggl(\frac1{\min (M_0^{3/4}, \sqrt{t})}+\epsilon\biggr).
\end{equation} 
Now we compute $\|v_R-I\|_{L^1(\Sigma_R\cap\overline{\calO}^c)}$.
We first focus on $\Sigma^{(3)}_{-1}\cap\overline{\calO}^c\cap\{ Im(z)>0\}$
When $1+\delta\le \frac{2t}{k}$, 
by \eqref{er92}, 
\begin{equation}
   \|v_R-I\|_{L^1(\Sigma^{(3)}_{-1}\cap\overline{\calO}^c\cap\{ Im(z)>0\})}
\le Ce^{-c_0(\epsilon)t}
\le \frac{C}{\sqrt{t}} \sqrt{1-\eta^2},
\end{equation} 
for large $t$ 
as $\sqrt{1-\eta^2}\ge C$ in this case.
When $1+\frac{M_0}{2^{1/3}k^{2/3}}\le \frac{2t}{k}<1+\delta$, 
from \eqref{er92}, 
\begin{equation}
\begin{split} 
  \|v_R-I\|_{L^1(\Sigma^{(3)}_{-1}\cap\overline{\calO}^c\cap\{ Im(z)>0\})}
&\le \int_{\Sigma^{(3)}_{-1}\cap\overline{\calO}^c\cap\{ Im(z)>0\}}
Ce^{-ct\sqrt{1-\eta^2}|z-\xi|^2}|dz|\\
&\le C\int_{\varrho}^{\infty} e^{-ct\sqrt{1-\eta^2}r^2} dr\\
%&\le C\int_{(t\sqrt{1-\eta^2})^{1/2}\varrho}^\infty
%e^{-cs^2} \frac{ds}{(t\sqrt{1-\eta^2})^{1/2}}\\
&\le \frac{C}{(t\sqrt{1-\eta^2})^{1/2}} 
e^{-ct\sqrt{1-\eta^2}\varrho^2}.
\end{split} 
\end{equation} 
But since, for $1+\frac{M_0}{2^{1/3}k^{2/3}}\le \frac{2t}{k}<1+\delta$,
\begin{equation}
  \sqrt{1-\eta^2}\ge \sqrt{1-\frac{k}{2t}}
=\frac1{(2t)^{1/3}}\biggl(\frac{k}{2t}\biggr)^{1/6}
\sqrt{k^{2/3}\biggl(\frac{2t}{k}-1\biggr)}
\ge \frac{CM_0^{1/3}}{t^{1/3}},
\end{equation} 
we obtain 
\begin{equation}\label{er103} 
  \|v_R-I\|_{L^1(\Sigma^{(3)}_{-1}\cap\overline{\calO}^c\cap\{ Im(z)>0\})}
\le \frac{C}{M_0^{3/4}}
e^{-c\epsilon^2M_0^{3/2}}\sqrt{1-\eta^2}.
\end{equation} 
By a similar computation, we obtain the same estimate 
for the other parts of $\Sigma_R\cap\overline{\calO}^c$.
Thus if we take $\epsilon$ small, and then take $M_0$, $t$ large, 
we obtain by \eqref{ee5.92}, \eqref{ee5.93} and \eqref{er103}, 
\begin{equation}\label{ee5.97}
   \|m^{(3)}_{11}-1\|_{L^1(\Sigma_R)} \le \alpha\sqrt{1-\eta^2}, 
\end{equation} 
with a constant $\alpha>0$ which can be taken to be arbitrarily small.
Therefore, from \eqref{ee5.46}, \eqref{ee5.92}, 
using \eqref{ee5.97}, 
we obtain (note \eqref{ee5.47})
for large $t$, 
\begin{equation}
   \log m_{11}(0)= \log m^{(3)}_{11}(0)+(1-\frac{\theta_c}{\pi})\log(1-s^2)
\le \alpha \sqrt{1-\eta^2} -c \sqrt{1-\eta^2}
\le -C\sqrt{1-\eta^2}, 
\end{equation} 
for some $C>0$, which is \eqref{ee1}.

\section{Multi-Painlev\`e Functions}\label{sec-ptons}

In this section we will show that the multi-interval case 
considered in Theorem \ref{thmmulti} is related to new classes of 
``multi-Painlev\'e function''.
As we will see, these functions describe the interaction 
of solutions of Painlev\'e equations in a way which is strongly 
reminiscent of the interaction of classical solitons.
We suggest the name ``Painlev\'etons''
or simply ``P-tons'' for these functions.
In this section we only illustrate a few of the properties of 
P-tons.
The general theory will be developed in a subsequent paper together 
with Alexander Its.

From Theorem \ref{thmmulti}, in the $k$ interval case, 
\begin{equation}
  \sum_{j=1}^k (s_j-s_{j-1}) \Kk_{n_j}(z,w)
= \frac{\sum_{l=0}^k f_l(z)g_l(w)}{z-w},
\end{equation} 
where 
\begin{align}
  &f=(f_0,\cdots,f_k)^T
= (s_k, (s_1-s_0)\varphi(z)z^{n_1}, \cdots, 
(s_k-s_{k-1})\varphi(z)z^{n_k})^T,\\
  &g=(g_0,\cdots,g_k)^T
= (2\pi i)^{-1}(1, -(\varphi(z)z^{n_1})^{-1}, \cdots, 
-(\varphi(z)z^{n_k})^{-1})^T.
\end{align} 
Thus by the integrable operator theory \cite{IIKS, deiftint}, 
the associated jump matrix $v$ on $\Sigma=\{|z|=1\}$ has the form 
\begin{equation}
  v=I-2\pi i fg^T= \begin{pmatrix}
1-s_k & s_k(\varphi z^{n_1})^{-1}& \hdots\hdots\hdots & 
s_k(\varphi z^{n_k})^{-1}\\
-(s_1-s_0)\varphi z^{n_1}  \\
\vdots & & \biggl(\delta_{pq}+(s_p-s_{p-1})z^{n_p-n_q}\biggr)_{1\le p,q\le k} \\
-(s_k-s_{k-1})\varphi z^{n_k}
\end{pmatrix} 
\end{equation} 

For purposes of illustration, we will only consider the case when $k=2$, 
\begin{equation}
  v=v^{(3)}= \begin{pmatrix} 
1-s_2 & s_2(\varphi z^{n_1})^{-1} & s_2(\varphi z^{n_2})^{-1} \\
-s_1\varphi z^{n_1} & 1+s_1 & s_1 z^{n_1-n_2} \\
-(s_2-s_1) \varphi z^{n_2} & (s_2-s_1)z^{n_2-n_1} & 1+ s_2-s_1
\end{pmatrix} 
\end{equation} 
and $\varphi=e^{t(z-z^{-1})}$ as in Introduction.
Observe now that when $s_1=0$, the jump matrix takes the form 
\begin{equation}
  v=v^{(3)}= \begin{pmatrix} 
1-s_2 & s_2(\varphi z^{n_1})^{-1} & s_2(\varphi z^{n_2})^{-1} \\
0&1&0\\
-s_2 \varphi z^{n_2} & s_2 z^{n_2-n_1} & 1+ s_2
\end{pmatrix}. 
\end{equation} 
Let $m^{(3)}$ be the solution of the $3\times 3$ RHP 
\begin{equation}
\begin{cases} 
  m^{(3)}_+=m^{(3)}_- v^{(3)}, \qquad z\in\Sigma,\\
  m^{(3)}\to I \qquad \text{as $z\to\infty$.}
\end{cases}
\end{equation} 
But it is clear that the $2\times 2$ matrix $m^{(2)}$ constructed from 
$m^{(3)}$ as follows, 
\begin{equation}
   m^{(2)}= \begin{pmatrix} 
 m^{(3)}_{11} & m^{(3)}_{13} \\
 m^{(3)}_{31} & m^{(3)}_{33}
\end{pmatrix} 
\end{equation} 
solves the RHP 
\begin{equation}
\begin{cases} 
  m^{(2)}_+=m^{(2)}_- \begin{pmatrix} 
1-s_2 & s_2(\varphi z^{n_2})^{-1} \\
-s_2 \varphi z^{n_2} & 1+s_2
\end{pmatrix}, \qquad z\in\Sigma,\\
  m^{(2)}\to I \qquad \text{as $z\to\infty$,}
\end{cases}
\end{equation}
which is an RHP which is algebraically equivalent to 
the RHP for Painlev\'e III (PIII)
which occurred in \cite{BDJ} : 
set 
\begin{equation}
\begin{cases} 
  \widetilde{m}^{(2)} = \begin{pmatrix} 
\sqrt{1+s_2}&0\\0&\frac1{\sqrt{1+s_2}} \end{pmatrix} 
m^{(2)}  \begin{pmatrix}
\sqrt{1+s_2}&0\\0&\frac1{\sqrt{1+s_2}} \end{pmatrix}
\qquad &|z|<1, \\
  \widetilde{m}^{(2)} = \begin{pmatrix}
\sqrt{1+s_2}&0\\0&\frac1{\sqrt{1+s_2}} \end{pmatrix}
m^{(2)}  \begin{pmatrix}
\frac1{\sqrt{1+s_2}}&0\\0&\sqrt{1+s_2} \end{pmatrix}
\qquad &|z|>1. 
\end{cases} 
\end{equation} 
Then $\widetilde{m}^{(2)}$ solves the RHP 
\begin{equation}
\begin{cases}
  \widetilde{m}^{(2)}_+=\widetilde{m}^{(2)}_- \begin{pmatrix}
1-s_2^2 & s_2(\varphi z^{n_2})^{-1} \\
-s_2 \varphi z^{n_2} & 1
\end{pmatrix}, \qquad z\in\Sigma,\\
  m^{(2)}\to I \qquad \text{as $z\to\infty$,}
\end{cases}
\end{equation}    
which is the RHP for PIII considered in \cite{FMZ}.
On the other hand, if $s_1=s_2=s$, then 
\begin{equation}
  v=v^{(3)}= \begin{pmatrix} 
1-s & s(\varphi z^{n_1})^{-1} & s(\varphi z^{n_2})^{-1} \\
-s\varphi z^{n_1} & 1+s & s z^{n_1-n_2} \\
0&0&1
\end{pmatrix} .
\end{equation} 
Now
\begin{equation}
   m^{(2)}= \begin{pmatrix}
 m^{(3)}_{11} & m^{(3)}_{12} \\
 m^{(3)}_{21} & m^{(3)}_{22}
\end{pmatrix}
\end{equation}
solves the RHP 
\begin{equation}
\begin{cases}
  m^{(2)}_+=m^{(2)}_- \begin{pmatrix}
1-s & s(\varphi z^{n_2})^{-1} \\
-s \varphi z^{n_2} & 1+s
\end{pmatrix}, \qquad z\in\Sigma,\\
  m^{(2)}\to I \qquad \text{as $z\to\infty$,}
\end{cases}
\end{equation}
which again is the (equivalent) RHP for PIII.
Also if we set $n_1=n_2=n$, 
\begin{equation}
  v=v^{(3)}= \begin{pmatrix} 
1-s_2 & s_2(\varphi z^{n})^{-1} & s_2(\varphi z^{n})^{-1} \\
-s_1\varphi z^{n} & 1+s_1 & s_1 \\
-(s_2-s_1) \varphi z^{n} & (s_2-s_1) & 1+ s_2-s_1
\end{pmatrix} .
\end{equation} 
Conjugating the solution $m^{(3)}$ of the RHP associated with $v^{(3)}$ by
\begin{equation}
m^{(3)}\mapsto \widetilde{m}^{(3)}=
\begin{pmatrix} 1&0&0\\ 0&1&0\\ 0&1&1 \end{pmatrix} 
m^{(3)} 
\begin{pmatrix} 1&0&0\\ 0&1&0\\ 0&1&1 \end{pmatrix}^{-1}
\end{equation} 
we find that $\widetilde{m}^{(3)}\to I$ as $z\to\infty$, and 
$\widetilde{m}^{(3)}$ solves a RHP with jump matrix 
\begin{equation}
  \widetilde{v}^{(3)}= \begin{pmatrix} 
1-s_2 & 0 & s_2(\varphi z^{n})^{-1} \\
-s_1\varphi z^{n} & 1 & s_1 \\
-s_2 \varphi z^{n} & 0 & 1+ s_2
\end{pmatrix}. 
\end{equation} 
It follows that necessarily
\begin{equation}
 (\widetilde{m}^{(3)}_{12} \quad \widetilde{m}^{(3)}_{22} 
\quad \widetilde{m}^{(3)}_{32} )^T
= (0 \quad 1 \quad 0)^T
\end{equation} 
and hence 
\begin{equation}
   m^{(2)}= \begin{pmatrix}
 \widetilde{m}^{(3)}_{11} & \widetilde{m}^{(3)}_{13} \\
 \widetilde{m}^{(3)}_{31} & \widetilde{m}^{(3)}_{33}
\end{pmatrix}
\end{equation}
solves the RHP 
\begin{equation}
\begin{cases}
  m^{(2)}_+=m^{(2)}_- \begin{pmatrix}
1-s_2 & s_2(\varphi z^{n})^{-1} \\
-s_2 \varphi z^{n} & 1+s_2
\end{pmatrix}, \qquad z\in\Sigma,\\
  m^{(2)}\to I \qquad \text{as $z\to\infty$,}
\end{cases}
\end{equation}
which is again the (equivalent) RHP for PIII.

The analogy with solitons is particularly clear if we consider $v^{(3)}$
in the edge scaling limit, 
\begin{equation}
  n_j=2\tau+t_j\tau^{1/3}, \qquad j=1,2 \quad ; \quad t_1<t_2,
\end{equation} 
as $\tau\to\infty$.
Then 
\begin{equation}
  v^{(3)}\bigl(-1+\frac{2iu}{\tau^{1/3}}\bigr) \mapsto 
\widehat{v}^{(3)}(u)= \begin{pmatrix} 
1-s_2 & s_2e^{2i\theta_1} & s_2e^{2i\theta_2} \\
-s_1e^{-2i\theta_1} & 1+s_1 & s_1e^{-2i(\theta_1-\theta_2)}  \\
-(s_2-s_1)e^{-2i\theta_2}  & (s_2-s_1)e^{2i(\theta_1-\theta_2)} & 1+ s_2-s_1
\end{pmatrix} 
\end{equation} 
on the real line, 
where 
\begin{equation}
 \theta_j=\frac43u^3+t_ju, \qquad j=1,2.
\end{equation} 
In addition to varying $s_1$, $s_2$, we can now vary $t_1$, $t_2$.
In particular, we can follow the trajectory of the 
solution of the RHP as $t_2$ moves from $t_1$ to $\infty$.
As $t_2\to t_1$, the solution becomes Painlev\'e II (PII) 
and as $t_2\to \infty$, it gives to another solution of PII, 
but now with a phase shift (see \cite{BDI}).
It is this behavior of P-tons, in particular, that is reminiscent of 
soliton interactions.

\section{Colored permutations}\label{sec-color}

     First, the definition:

Let $\pi$ be an $m$-colored permutation 
(see, e.g., \cite{Rains}), and assume the colors are
indexed by $0,1,\dots m-1$.  Let $S$ be a subsequence of length $l$ of
$\pi$ which is a union of monochromatic increasing subsequences;
let $k_i$ be the number of these sequences having color $i$, and set
$k=\sum_i k_i$.  Note that the monochromatic increasing subsequences may
be empty, but the color of empty subsequences still matters.  We assign to
$S$ the following score:
\begin{equation}
ml+{k+1\choose 2}+\sum_{0\le i\le m-1} (ik_i-m{k_i+1\choose 2})
\end{equation}
Now, let $l_k(\pi)$ be the maximum score over all unions of $k$ monochromatic
increasing subsequences (note $l_0(\pi)=0$).  We then define
\begin{equation}
\lambda_k(\pi):=l_k(\pi)-l_{k-1}(\pi).
\end{equation}

\begin{lem}
Let $\lambda^{(i)}_k(\pi)$ be the partition associated to just the
$i$-colored subsequence of $\pi$.  Then $\lambda_k(\pi)-k$
is simply the $k$th largest of the numbers $m(\lambda^{(i)}_j(\pi)-j)+i$.
Moreover, if $\pi$ has length $n$, then $\lambda_k(\pi)$ is a partition
of $mn$.
\end{lem}

\begin{proof}
Fix a composition $k_i$, and consider the largest score associated to that
composition.  Clearly, we can maximize the score for each color independently;
we thus obtain:
\begin{equation}
ml + {k+1\choose 2} + \sum_{0\le i\le m-1} (ik_i-m{k_i+1\choose 2})
=
{k+1\choose 2} +
\sum_{0\le i\le m-1} \sum_{1\le j\le k_i} (m(\lambda^{(i)}_j(\pi)-j)+i).
\end{equation}
Now, for a fixed value of $k$, this is clearly maximized when
the values $m(\lambda^{(i)}_j(\pi)-j)+i$ occurring in the sum are
chosen to be as large as possible.  Plugging the resulting value of $l_k(\pi)$
into the formula for $\lambda_k(\pi)$, we obtain the first claim.

Note that the numbers $m(\lambda^{(i)}_j(\pi)-j)+i$ are all different
(the congruence class modulo $m$ depends on the color, and the numbers are
distinct within a given color).  Furthermore, we readily verify that for
each congruence class, the number of negative numbers not occurring in the
set is equal to the number of nonnegative numbers occurring in the set.  We
thus conclude that $\lambda_k(\pi)$ is indeed a partition.  It remains
to verify that $\sum_k \lambda_k(\pi)=mn$; in other words, $l_k(\pi)=mn$
for $k$ sufficiently large.  Choose $k$ such that $\pi$ is a union of $k$
increasing subsequences, and consider $l_{mk}(\pi)$.  We readily verify
that the term
\begin{equation}
\sum_{0\le i\le m-1} (ik_i-m{k_i+1\choose 2})
\end{equation}
is maximized when all $k_i$ are equal to $k$, and thus the optimal score
differs from $mn$ by
\begin{equation}
{mk+1\choose 2} + \sum_{0\le i\le m-1} (ik-m{k+1\choose 2})
=
0.
\end{equation}
\end{proof}

\begin{rem}
An alternate approach is to define $\lambda_k(\pi)$ via the
Schensted correspondence for rim-hook permutations given in
\cite{StantonWhite}, at which point the lemma follows immediately.
The fact that the rim-hook correspondence splits into $m$ ordinary
correspondences gives the increasing subsequence interpretation above.
\end{rem}

     Now, suppose we choose $n$ randomly according to a Poisson law of
mean $mt^2$, and then choose an $m$-colored permutation of length $n$ at
random.  Equivalently, take $m$ independent Poisson processes in the
unit square (one for each color), and convert the resulting point set
to a colored permutation.  We thus see that the resulting random partitions
$\lambda^{(i)}_j(\pi)$ are independent, and are all distributed according to
the law for ordinary permutations.  In particular, we obtain the following
correlation kernel:
\begin{equation}
\Ss^{(m)}(a,b)
=
\sum_{k\ge 1} (\varphi^{-1})_{(a+k)/m}\varphi_{(b+k)/m}
\end{equation}
where
\begin{equation}
\varphi(z)= e^{t(z-z^{-1})}.
\end{equation}
(Recall from Corollary \ref{cor-id} that $\varphi_a$ and $(\varphi^{-1})_a$
are $0$ for $a$ non-integral.)

Now by using Corollary \ref{cor-id} and Theorem \ref{thmmom} 
for the convergence of moments for the ordinary permutations, 
we obtain the convergence of moments for $\lambda_k$'s 
in the colored permutation setting.
More precisely, as in \eqref{ec3.5}, 
there is a limiting distribution $F^{color(m)}$ such that 
\begin{equation}
  \lim_{N\to\infty}  \Exp^{color(m)}_N
\biggl( \prod_{j=1}^k\biggl(
\frac{\lambda_j-2\sqrt{mN}}{m^{2/3}(mN)^{1/6}}\biggr)^{a_j}
\biggr) = 
\Exp^{color(m)}\bigl(x_1^{a_1}\cdots x_k^{a_k}\bigr)
\end{equation} 
where $\Exp^{color(m)}_N$ denotes the expectation with respect to the natural 
counting measure on the colored permutations (see \cite{Rains}), 
and $\Exp^{color(m)}$ is the expectation with respect to $F^{color(m)}$.
%the limiting 
%distribution $F^{color(m)}$ in \eqref{eq2.4}, \eqref{ei5}.
The function $F^{color(m)}(x_1,\cdots,x_k)$ has the following meaning in 
terms of GUE. Take $m$ random GUE matrices of size $N$ at random, then 
superimpose their eigenvalues. We denote the largest of those superimposed 
numbers by $z_1(N)$, the second largest by $z_2(N)$, and so on.
Then $F^{color(m)}(x_1, \cdots, x_k)$ is the limiting distribution of 
$z_1, \cdots, z_k$ as $N\to\infty$, after appropriate 
centering and scaling.

A number of other statistical systems which are currently of interest 
can also be analyzed by the methods of this paper.
In particular, we have in mind the random word problem 
\cite{TW:randomwords, kurtj:disc, ITW1, ITW2}, 
certain 2-dimensional growth models 
\cite{kurtj:shape}, and also the so-called ``digital boiling model''
\cite{GTW}.

For example, in the growth model considered by Johansson 
in \cite{kurtj:shape}, 
let $\sigma=\cup_{j=1}^k \sigma_j$ be a union of $k$ disjoint increasing 
paths $\sigma_j$ in the model.
Let $L^{(k)}(\sigma)$ be the sum of the lengths of the paths $\sigma_j$, 
and let $L^{(k)}=\max_{\sigma} L^{(k)}(\sigma)$. 
We define $\lambda_k=L^{(k)}-L^{(k-1)}$.
The joint probability distribution for $\lambda_1, \cdots ,\lambda_k$ 
can be obtained \cite{kurtj:shape} by various differentiations 
of $\det(1+\sum_{j=1}^k s_j \chi_{[n_j.n_{j-1})}\Ss)$ 
with respect to $s_1,\cdots,s_k$ as in \eqref{ei10} with 
$\varphi$ now given by $\varphi(z)=(1+\sqrt{q}z)^M(1+\sqrt{q}z^{-1})^{-N}$.
But now by Theorem \ref{thmmulti}, 
$\det(1+\sum_{j=1}^k s_j \chi_{[n_j.n_{j-1})}\Ss)$ can be expressed in terms
of the determinant of an integrable operator as in \eqref{ethmmulti}.
This opens up the possibility for the asymptotic analysis of the convergence 
of moments for the joint distribution.
However, the associated RHP has a new feature, namely 
the weight function is non-real,  
which has not yet been addressed in general (however, see \cite{KMM}).
There are similar formulae for random words and digital boiling.

\bibliographystyle{plain}
\bibliography{paper9}

\end{document}